\documentclass[journal]{IEEEtran}
\usepackage{etex}
\reserveinserts{100}
\usepackage{mathtools}
\usepackage[dvipsnames]{xcolor}
\usepackage{graphics, theorem, times, amsfonts, graphicx, amssymb, cite}
\usepackage[outdir=./]{epstopdf}
\usepackage{tikz}
\usetikzlibrary{shapes,arrows}
\usepackage{pgfplots}
\usepackage{color}
\usepackage{setspace}
\usepackage{hyperref}
\usepackage{multirow}
\usepackage{rotating}
\usepackage{comment}
\usepackage{flushend}
\usepackage[font=small]{caption}
\usepackage{subcaption}
\usepackage[affil-it]{authblk}
\usepackage{bm}
\usepackage{algorithm}
\usepackage{algorithmic}
\usepackage{enumerate}
\usepackage{morefloats}
\usepackage{setspace}

\input{mysymbol.sty}
\usepackage{needspace}

\newcommand{\QED}{\hfill\ensuremath{\blacksquare}}

\def\whp{\text{w.h.p}}
\def\E{\mathbb{E}}

\newtheorem{assumption}{\hspace{0pt}\bf Assumption}
\newtheorem{lemma}{\hspace{0pt}\bf Lemma}
\newtheorem{proposition}{\hspace{0pt}\bf Proposition}
\newtheorem{example}{\hspace{0pt}\bf Example}

\newtheorem{theorem}{\hspace{0pt}\bf Theorem}
\newtheorem{corollary}{\hspace{0pt}\bf Corollary}

\newtheorem{remark}{\hspace{0pt}\bf Remark}

\title{Learning in Wireless Control Systems over Non-Stationary Channels}
\author{Mark Eisen, Konstantinos Gatsis, George J. Pappas, and Alejandro Ribeiro
\thanks{{Supported by ARL DCIST CRA W911NF-17-2-0181 and Intel Science and Technology Center for Wireless Autonomous Systems. The authors are with the Department of Electrical and Systems Engineering, University of Pennsylvania. Email at: \{maeisen, kgatsis, pappasg, aribeiro\}@seas.upenn.edu.{This paper expands upon the results and presents convergence proofs that are referenced in \cite{ICASSPpaper,ACCpaper2}}.
}}}

\begin{document}
\thispagestyle{empty}
\maketitle

\begin{abstract}
This paper considers a set of multiple independent control systems that are each connected over a non-stationary wireless channel. The goal is to maximize control performance over all the systems through the allocation of transmitting power within a fixed budget. This can be formulated as a constrained optimization problem examined using Lagrangian duality. By taking samples of the unknown wireless channel at every time instance, the resulting problem takes on the form of empirical risk minimization, a well-studied problem in machine learning. Due to the non-stationarity of wireless channels, optimal allocations must be continuously learned and updated as the channel evolves. The quadratic convergence property of Newton's method motivates its use in learning approximately optimal power allocation policies over the sampled dual function as the channel evolves over time. Conditions are established under which Newton's method learns approximate solutions with a single update, and the subsequent sub-optimality of the control problem is further characterized. Numerical simulations illustrate the near-optimal performance of the method and resulting stability on a wireless control problem.
\end{abstract}


\section{Introduction}

The recent developments in autonomy in industrial control environments, teams of robotic vehicles, and the Internet-of-Things have motivated intelligent design of wireless systems.  Even though wireless communication facilitates connectivity, it also introduces uncertainty that may affect stability and performance. To guarantee performance and safety of the control application it is common to employ model-based approaches. However wireless communication is naturally uncertain and time-varying due to effects that are not always amenable to modeling, such as mobility in the environment. In this paper we propose an alternative learning-based approach, where autonomy relies on collected channel samples to optimize performance in a non-stationary environment. The connection between the two approaches is based on the observation that a sampled version of the model-based design approach can be cast as an empirical risk minimization (ERM) problem, a typical machine learning problem. Even so, standard techniques developed for solving ERM problems in machine learning do not address the additional challenges present in wireless autonomous systems, namely the non-stationarity of sample distributions.

The traditional model-based approach is motivated by the desire to build wireless control systems with stability and optimal performance. To counteract channel uncertainties it is natural to include a model of the wireless communication, for example an i.i.d. or Markov link quality, alongside the model of the physical system to be controlled. These models have been valuable to help analysis and control/communication design. For example, one can characterize that it is impossible to estimate and/or stabilize an unstable plant if its growth rate is larger than the rate at which the link drops packets~\cite{Sinopoli_intermittent, Schenato_foundations, hespanha2007survey, gupta2007optimal}, or below a certain channel capacity~\cite{Tatikonda, Sahai_anytime}. Additionally models facilitate the design of controllers~\cite{elia2005remote,imer2006optimal,cloosterman2010controller}, as well as the allocation of communication resources to optimize control performance, for example power allocation over fading channels with known distributions~\cite{GatsisEtal14, Quevedo_Kalman}, or event-triggered control~\cite{Rabi_scheduling, mazo2011decentralized, Event_triggered_intro, araujo2014system, mamduhi2014event}.

In practice wireless autonomous systems operate under unpredictable channel conditions following unknown time-varying distributions. While one approach would be to estimate the distributions using channel samples and then follow the above model-based design approach, in this paper we propose an alternative learning-based approach which bypasses the channel-modeling phase. We exploit channel samples taken from the time-varying channel distributions with the goal to learn directly the solution to communication design problems. To apply this approach we exploit a connection between the model-based and the learning-based design problems. Existing works \cite{GatsisEtal15, GatsisEtal18, chen2017stochastic} study related problems in multiple-access wireless control systems and resource allocation problems in wireless systems but under a stationary channel distribution. These works generally employ first-order stochastic methods, which have slow convergence rates and hence not suitable for the present framework. A significant challenge remains in how to continuously learn optimal policies over a wireless channel that is time-varying. This shortcoming of existing sample-based approaches used in \cite{GatsisEtal15, GatsisEtal18,chen2017stochastic} and more general machine learning scenarios motivates the higher-order learning approach proposed in this paper. Some existing machine learning methods account for nonstationarity by optimizing an averaged objective over all time \cite{mokkadem2011generalization,kingma2014adam,besbes2015non}. Our approach differs in that we seek and track optimality locally with respect to the current channel distribution at every time epoch. 

In this paper we consider a wireless autonomous system where the design goal is to maximize a level of control performance for multiple systems while meeting a desired transmit power budget over the wireless channel (Section~\ref{sec:formulation}). The wireless channel is modeled as a fading channel with a time-varying and unknown distribution, and only available through samples taken over time. We derive in Section \ref{sec:formulation2} a wireless control problem that finds optimal power allocation policies for an individual time epoch where the wireless channel distribution does not change, and then proceed to derive the Lagrange dual (Section \ref{sec:dual_form}). We show in Section~\ref{sec:ERM} that the dual of the power allocation problem can be rewritten using channel samples as an empirical risk minimization problem, a common machine learning problem in which an expected loss function over an unknown distribution is approximated by optimized over a set of samples. Here the risk is loosely related to how far the current solution is from the desired optimal power allocation. 

Because the wireless channel is varying over time, we develop a new approach to solving a sequence of ERM problems. We collect and store a window of channel samples taken from consecutive distributions to reduce sampling complexity and employ Newton's method to learn new policies quickly (Section~\ref{sec_nonstat}). More specifically, the quadratic convergence rate of Newton's method is shown to be sufficient to find approximate solutions to slowly varying objectives with a single update. Using Newton's method, we propose an algorithm that uses channel samples to approximate the solution of a power allocation wireless control problem over a non-stationary channel. We prove that, under specific conditions, the algorithm reaches an approximately optimal point in a single iteration of Newton's method (Section \ref{sec_convergence}). This result establishes both a suboptimality bound with respect to the sampled problem (Section \ref{sec_erm_conv}) as well as with respect to control performance metric in the wireless control problem (Section \ref{sec_subopt}). We additionally show a stability result for a particular problem description common in wireless control systems (Section \ref{sec_stability}) and provide considerations for practical implementation of the method (Section \ref{sec_implementation}). These results are further demonstrated in a numerical demonstration of learning power allocation policies across multiple control systems over a time-varying channel (Section \ref{sec_numerical_results}).

\section{Wireless Control Problem}\label{sec:formulation}

\begin{figure}
\centering
\input{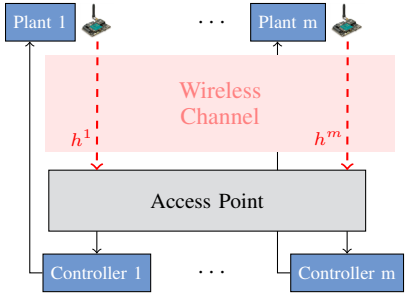}
\caption{Wireless control system. Plants communicate state information to access point/controllers over wireless medium.}
\label{fig_channel}
\end{figure}

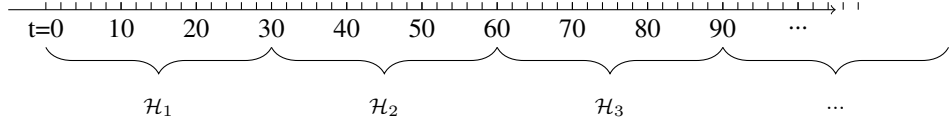
\begin{figure*}
\centering
\usetikzlibrary{decorations.pathreplacing,angles,quotes,calc}

\begin{tikzpicture}
\draw (0 ,0) node[anchor=north] {t=0} -- (0,0.1);
\foreach \x in {1,2,...,9} { 
  
  \draw (\x ,0) node[anchor=north] {\x0} -- (\x,0.1);
};
\foreach \x in {0.2,0.4,...,10.8} { 
  
  \draw (\x ,0) node[anchor=north] {} -- (\x,0.1);
};
\draw (10 ,0) node[anchor=north,yshift=-.1cm] {...} -- (10,0.1);

\draw[->] (-0.5,0) -- (10.5,0);
\draw [decorate,decoration={brace,amplitude=10pt,mirror}]
(0,-0.5) -- (3,-0.5) node [black,midway,xshift=-0,yshift=-0.8cm] 
{\footnotesize $\mathcal{H}_1$};

\draw [decorate,decoration={brace,amplitude=10pt,mirror}]
(3,-0.5) -- (6,-0.5) node [black,midway,xshift=-0,yshift=-0.8cm] 
{\footnotesize $\mathcal{H}_2$};

\draw [decorate,decoration={brace,amplitude=10pt,mirror}]
(6,-0.5) -- (9,-0.5) node [black,midway,xshift=-0,yshift=-0.8cm] 
{\footnotesize $\mathcal{H}_3$};

\draw [decorate,decoration={brace,amplitude=10pt,mirror}]
(9,-0.5) -- (12,-0.5) node [black,midway,xshift=-0,yshift=-0.8cm] 
{\footnotesize $...$};
\end{tikzpicture}
\caption{Time axis showing evolution of time $t$ and epochs $k$. Each channel distribution $\ccalH_k$ is stationary for a set of time instances.}
\label{fig_time_scale}
\end{figure*}

We consider a wireless control problem (WCP) with $m$ independent control systems labeled $i=1,\hdots,m$, as shown in Figure \ref{fig_channel}. Each control system/agent $i$ communicates at time $t$ its state $x^i_t$ over a wireless channel in order to close a loop and maximize a level of control performance. In particular, system $i$ tries to close the control loop over the wireless channel by transmitting with power level $p^{i} \in  [0,p_0]$. Due to propagation effects the channel fading conditions that each system $i$ experiences, denoted by $h^{i} \in \reals_+$, change unpredictably over time~\cite[Ch.~3]{Goldsmith_Wireless_Communications}. Together, the channel fading $h^i$ and transmit power $p^i$ determine the signal-to-noise ratio (SNR) at the receiver for system $i$, which in turn affects the probability of successful decoding of the transmitted packet at the receiver. We consider a function $q( h^{i}, p^{i})$ that, given a current channel state and transmit power, determines the probability of successful transmission and decoding of the transmitted packet -- see, e.g.,~\cite{Quevedo_Kalman,GatsisEtal14} for more details on this model. Transmission are assumed on different frequencies/bands and are not subject to contention -- see~\cite{GatsisEtal15,GatsisEtal18} for alternative formulations. 

 Because these fading conditions vary quickly and unpredictably, they can be modeled as independent random variables drawn from distribution $\ccalH$ that itself is non-stationary, or time-varying. Channel fading is assumed constant during each transmission slot and it is independently distributed over time slots (block fading). Furthermore, the channel distribution $\ccalH$ may vary across \emph{time epochs}, but will in general be stationary within a single time epoch.  In particular, consider an epoch index $k=0,1,\hdots$ that specifies a particular channel distribution $\ccalH_k$ with realization $h^i_k$ for system $i$. In Figure \ref{fig_time_scale}, we display a time axis rendering of this model. The state variables change at each transmission slot $t$, while the channel changes at scale $k$, which will in general contain multiple time steps. This is to say that we assume that the channel distribution $\ccalH_k$ changes at a rate slower than the system evolution, and that within a single time epoch the channel is effectively stationary.
 
 We proceed to derive a formal description of the wireless control problem of interest within a single time epoch, where the channel is assumed stationary. In Section \ref{sec_nonstat} we extend this formulation to the non-stationary setting.

\subsection{WCP in single epoch}\label{sec:formulation2}
Within a particular time epoch $k$ with channel distribution $\ccalH_k$, we can derive a formulation that characterizes the optimal power allocations between the $m$ control systems so as to maximize the aggregate control performance across all systems, where $p_0$ reflects a maximum transmission power of the system. Given a random channel state $h_k^i \in \reals_+$ drawn from the distribution $\ccalH_k$. We wish to determine the amount of transmit power $p_k^i(h_k^i): \reals_+ \rightarrow [0,p_0]$ to be used when attempting to close its loop---see \cite{GatsisEtal14} for details. We note that we are looking for transmit power as a function of current channel conditions, as the power necessary to close the loop will indeed change with channel conditions.  We assume the current channel gain $h_k^i$ is available at the transmitter at each slot, as this can generally be obtained via short pilot signals---see \cite{GatsisEtal14}. Then the probability of closing the loop is given by the value
\begin{equation}\label{eq:y_definition}
y_k^i := \E_{h_k^i} \left\{q(h_k^i,p_k^i(h_k^i))\right\}.
\end{equation}
The variable $y_k^i \in [0,1]$ is the expectation of successful transmission over the channel distribution $\ccalH_k$. 

Using the variable $y_k^i$ we use a monotonically increasing concave function $J^i: [0,1] \rightarrow \reals$ that returns a measure of control system performance as a function of the probability of successful transmission. Such a function can take on many forms and, in general, can be derived  in relation to the particular control task of interest.  In the following example, we derive such a measure for a typical wireless control problem setting, namely the quadratic control performance of a switched linear dynamical system -- see, e.g.,~\cite{Schenato_foundations,gupta2007optimal}.

\begin{example}\label{example1}
	Consider for example that a control system $i$ is a scalar linear dynamical system of the form	
	\begin{equation}\label{eq:system}
	x^i_{t+1} = A^i_o x^i_t + B^i u^i_t + w^i_t
	\end{equation}
	where $x^i_t \in \reals$ is the state of the system at transmission time $t$, $A^i_o$ is the open loop (potentially unstable) dynamics of the system, $u^i_t \in \reals$ is the control input applied to the system at time $t$, and $w^i_t$ is some zero-mean i.i.d. disturbance process with variance $W^i$. Consider a given linear state feedback is applied to the system as the control input when a transmission is successful, i.e., 
	\begin{equation}\label{eq:system}
	u^i_{t} = \left\{ \begin{array}{ll} K^i x^i_t & \text{if loop closes} \\ 0 & \text{otherwise} \end{array} \right.
	\end{equation}
	As a result, the system switched between an open loop mode $A^i_o$ and a closed loop stable mode $A^i_c = A^i_o + B^i K^i$, as in 
	\begin{equation}\label{eq:system}
	x^i_{t+1} = \left\{ \begin{array}{ll} A^i_o x^i_t + w^i_t & \text{if loop closes} \\ A^i_o x^i_t + w^i_t & \text{otherwise} \end{array} \right.
	\end{equation}
	The goal is to regulate the system state close to zero, i.e., the system attempts to close the loop at a high rate in order to minimize an expected quadratic control cost objective of the form
	\begin{equation}\label{eq:quadratic_cost}
	\lim_{N \rightarrow \infty} \frac{1}{N} \sum_{t=0}^{N-1}\mathbb{E} (x^i_t)^2
	\end{equation}
	Assuming the control loop in \eqref{eq:system} is closed with the success probability $y_k^i$ in \eqref{eq:y_definition} \textit{at all time steps}, it is possible to express the above cost explicitly as a function of $y_k^i$. Using the system dynamics \eqref{eq:system}, the variance of the system state satisfies the recursive formula
	\begin{equation}\label{eq_system_d}
	\mathbb{E}(x^i_{t+1})^2 = y_k^i \, (A^i_c)^2 \, \mathbb{E}(x^i_t)^2 + (1-y_k^i) \, (A^i_o)^2 \, \mathbb{E}(x^i_t)^2 + W^i
	\end{equation} 
	that is, with probability $y_k^i$ the variance grows according to the open loop dynamics, and with probability $1-y_k^i$ the variance shrinks according to the closed loop stable dynamics. 
	
	Operating recursively and using the geometric series sum, we can rewrite the variance at time $t$ as
	\begin{align}
	\mathbb{E}(x^i_{t})^2 &= [y_k^i \, (A^i_c)^2 +\,(1-y_k^i) \, (A^i_o)^2]^{t} \mathbb{E}(x^i_0)^2 \\
	&+ W^i\frac{ 1- [y_k^i \, (A^i_c)^2 +\,(1-y_k^i) \, (A^i_o)^2]^{t} }{1-[y_k^i \, (A^i_c)^2 + \,(1-y_k^i) \, (A^i_o)^2] } .
	\end{align} 
	As follows from the above expression, the system is stable, i.e., the variance is bounded, if the packet success rate satisfies $[y_k^i \, (A^i_c)^2 +\,(1-y_k^i) \, (A^i_o)^2] <1$ so that the sum above is bounded -- see also~\cite{Schenato_foundations,gupta2007optimal}. In that case, the state variance as well as the average \eqref{eq:quadratic_cost} converge to the same limit value, which we can define as our 
	control performance function
	\begin{equation}\label{eq_j_ex1}
	J^i(y_k^i) = -\frac{W^i}{1 - \left[y_k^i (A^i_c)^2 + (1-y_k^i) (A^i_o)^2\right]}
	\end{equation}
	This control performance function satisfies the assumption of concavity, and it is also monotonically increasing because we have added the negative sign in front of the expression. It is also possible to extend this analysis to include a cost on the control input, as is common in the Linear Quadratic Control problem, i.e., replace the cost in \eqref{eq:quadratic_cost} with $\mathbb{E} (x^i_t)^2 + (u^i_t)^2$.
\end{example}

\begin{remark}\label{remark_example}\normalfont
In Example \ref{example1}, observe that the control system performance in \eqref{eq:quadratic_cost} is a long term objective asymptotically for $t \rightarrow \infty$. As the channel fading distribution $\ccalH_k$ will change unpredictably in the future it is hard to define an accurate value of this control performance. As a surrogate, in the above example we write a control system performance in~\eqref{eq_j_ex1} with respect to the current channel distribution $\ccalH_k$, i.e., as if this channel distribution is stationary and will not change in the future. Later, in Section~\ref{sec_stability} we argue that this approximation and the power allocation algorithm we develop can indeed guarantee system stability.
\end{remark}


To derive the full formulation of the wireless control problem for current channel distribution $\ccalH_k$, we first define using boldface vectors the set of $m$ channel states $\bbh_k := [h_k^1; h_k^2; \hdots; h_k^m] \sim \ccalH_k^m$ observed by the control systems and the set of power allocation policies $\bbp_k(\bbh_k) := [p_k^1(h_k^1); p_k^2(h_k^2); \hdots; p_k^m(h_k^m)] \in \ccalP := [0,p_0]^m$. We further define the vector of transmission probabilities at specific channel states $\bbq(\bbh_k, \bbp_k(\bbh_k)) := [q(h_k^1,p_k^1(h_k^1)); \hdots; q(h_k^m,p_k^m(h_k^m))]$ and expected transmission probabilities $\bby_k := [y_k^1; y_k^2; \hdots; y_k^m]$ from \eqref{eq:y_definition}. The goal is to select $\bbp_k(\bbh_k)$  whose expected aggregate value is within a maximum power budget $p_{\max}$ while maximizing the total system performance $\sum J^i$ over $m$ agents.
Because $J^i$ is monotonically increasing, we can relax the equality in \eqref{eq:y_definition} to an inequality constraint and write the following optimization problem.
\begin{align}\tag{WCP$_k$}\label{eq_power_problem_y}
\{\bbp_k^*(h), \bby_k^*\} &:= \argmax_{\bbp_k \in \ccalP,\bby_k \in \reals^m} J(\bby_k) := \sum_{i=1}^m J^i (y_k^i) \\
&\st \qquad    \bby_k \leq \E_{\bbh_k} \left\{\bbq(\bbh_k,\bbp_k(\bbh_k))\right\},  \nonumber \\
&\quad \qquad \sum_{i=1}^m \E_{h_k^i} (p_k^i(h_k^i)) \leq p_{\max} \nonumber 
\end{align}

The problem in \eqref{eq_power_problem_y} states the optimal power allocation policy $\bbp_k^*(\bbh_k)$ is the one that maximizes the expected aggregate control performance over channel states while guaranteeing that the expected total transmitting power is below an available budget $p_{\max}$. We stress that this only provides the optimal policy with respect to a particular channel distribution $\ccalH_k$. In the non-stationary wireless setting we are interested in solving \eqref{eq_power_problem_y} for all $k$.

\subsection{Dual formulation of \eqref{eq_power_problem_y}}\label{sec:dual_form}
Solving this optimization problem directly has a number of significant challenges. The first is that the problem is non-convex, in particular due to the first constraint in \eqref{eq_power_problem_y}. The second challenge is that the problem is optimized over an infinite-dimensional variable $\bbp_k(\bbh_k)$. It is very difficult to solve such a problem if there is no assumed parameterization of $\bbp_k^*(\bbh_k)$. We can show, however, from a result in \cite{ribeiro2012optimal} that a naturally occurring parameterization of $\bbp_k^*(\bbh_k)$ indeed can be derived from Lagrangian duality theory.

We proceed then to derive the dual problem from the constrained problem in \eqref{eq_power_problem_y}. To simplify the presentation, we first introduce a set of augmented variables, denoted with tildes. Define the augmented vectors $\check{\bbq}(\bbh_k,\bbp_k(\bbh_k)) \in \reals^{m+1}$ and $\check{\bby}_k \in \reals^{m+1}$ as
\begin{align}\label{eq_augmented_vars}
\check{\bbq}(\bbh_k,\bbp_k(\bbh_k)) := \begin{bmatrix}
q(h_k^1,p_k^1(h_k^1)) \\ \vdots \\ q(h_k^m,p_k^m(h_k^m))\\ -\sum_{i=1}^mp_k^i(h_k^i)] \end{bmatrix} \quad
\check{\bby}_k := \begin{bmatrix}
y_k^1 \\ \vdots \\ y_k^m \\ -p_{\max} \end{bmatrix}.
\end{align}
The augmented $\check{\bbq}(\bbh_k,\bbp_k(\bbh_k))$ includes transmission probabilities augmented with the total power allocation while $\check{\bby}_k$ includes auxiliary variables augmented with total power budget. Using this new notation, the Lagrangian function is formed as 
\begin{align}\label{eq_lagrangian}
\ccalL_k(\bbp_k(\bbh_k),\bby_k,\bbmu_k) &:= \sum_{i=1}^m J^i (y_k^i) \\+ & \bbmu_k^T\left(\E_{\bbh_k} \check{\bbq}(\bbh_k,\bbp_k(\bbh_k)) - \check{\bby}_k\right), \nonumber
\end{align}
where $\bbmu_k  := [\mu_k^1; \hdots; \mu_k^m; \tdmu] \in \reals_{+}^{m+1}$ contains the dual variables associated with each of the $m+1$ constraints in \eqref{eq_power_problem_y}. From the Lagrangian function in \eqref{eq_lagrangian}, the Lagrangian dual loss function is defined as $L_k(\bbmu_k) := \max_{\bbp_k,\bby_k} \ccalL_k(\bbp_k(\bbh_k),\bby_k, \bbmu_k)$---see, e.g., \cite{boyd2004convex}---and the corresponding dual problem as
\begin{align}\label{eq_dual_problem}
\tbmu_k^* &:= \argmin_{\bbmu_k \geq 0} L_k(\bbmu_k) \\
L_k(\bbmu_k) &:=  \E_{\bbh_k}\! \!\!\left\{\max_{\bbp_k,\bby_k}  \sum_{i=1}^m J^i (y_k^i)\!\! + \bbmu_k^T\left( \check{\bbq}(\bbh_k,\bbp_k(\bbh_k)) \!\! - \check{\bby}_k\right)\right\}. \nonumber
\end{align}
Note in \eqref{eq_dual_problem} that the expectation operator and maximization were exchanged without loss of generality---see, e.g. \cite[Proposition 2]{GatsisEtal15}. It is important to stress here the connection between the dual problem in \eqref{eq_dual_problem} with the original problem in \eqref{eq_power_problem_y}. While \eqref{eq_power_problem_y} is indeed not convex, problems of this form can be shown to exhibit zero duality gap under the technical assumption that the primal problem is strictly feasible and that the channel probability distribution is non-atomic  \cite{ribeiro2012optimal}. This implies that the optimal primal variable $\bbp_k^*(\bbh_k)$ in \eqref{eq_power_problem_y} can be recovered from the optimal dual variable $\tbmu_k^*$ in \eqref{eq_dual_problem}. Thus, the power allocation policy for each agent $i$ is found indirectly by solving \eqref{eq_dual_problem} and recovering as 
\begin{align}\label{eq_primal_recovery}
p_k^i(h_k^i,\bbmu_k) &= \argmax_{p_k^i \in [0,p_0]}  \mu_k^i q(h_k^i,p_k^i(h_k^i)) -  \tdmu p_k^i(h_k^i) ,  \\
y_k^i(\bbmu_k) &= \argmax_{y_k^i}  J^i (y_k^i) - \mu_k^i y_k^i. \label{eq_primal_recovery_2}
\end{align}
 The optimal policy is subsequently recovered using the optimal dual variable as $\bbp_k^*(\bbh_k) := [p_k^1(h_k^1,\tbmu_k^*);\hdots;p_k^m(h_k^m,\tbmu_k^*)]$. Observe that the problem in \eqref{eq_dual_problem} is a simply constrained stochastic problem that is known to always be convex from duality theory, and can be solved efficiently with a variety of projected stochastic descent methods \cite{bottou2010large,duchi2011adaptive,defazio2014saga,GatsisEtal15, GatsisEtal18}. Thus, the non-convex, infinite-dimensional optimization problem in \eqref{eq_power_problem_y} can be solved indirectly but exactly with the convex, finite-dimensional problem in \eqref{eq_dual_problem}.

\begin{remark}\label{remark_problem}\normalfont
The problem formulation given in \eqref{eq_power_problem_y} that we use in this paper assumes there is a fixed power budget and the metric to be optimized is a measure of control performance. An alternative formulation of resource allocation that may be more relevant in some settings would instead fix a bound on the required control performance, typically derived from a stability margin for the control system. Here the objective would instead be to minimize total power usage, subject to the constraint on control performance. Indeed, these two problems are very similar when reformulated in the dual domain, and can thus be studied almost identically as such. We specifically focus on the problem in \eqref{eq_power_problem_y} in this paper but stress that all the results will apply to this alternative problem as well.
\end{remark}

\section{ERM Formulation of \eqref{eq_power_problem_y}}\label{sec:ERM}

The stochastic program in \eqref{eq_dual_problem} features an objective that is the expectation taken over a random variable, and can thus be  considered as a particular case  of the empirical risk minimization (ERM) problem. Empirical risk minimization is a common problem studied in machine learning due to its ubiquity in training classifiers,  and the same structure appears naturally in the dual formulation of the WCP. A generic ERM problem considers a convex loss function $f(\bbmu_k, \bbh_k)$ of a decision variable $\bbmu_k \in \reals^{m+1}$ and random variable $\bbh_k$ drawn from distribution $\ccalH_k$ and seeks to minimize the expected loss $L_k(\bbmu_k) := \mathbb{E}_{\bbh_k} [ f(\bbmu_k,\bbh_k)]$. For the WCP in \eqref{eq_power_problem_y}, we rewrite the loss function $L$ and associated ERM problem in terms of a function $f(\bbmu_k,\bbh_k)$ using its dual as
\begin{align}\label{eq_erm_problem}
&\tbmu_k^* := \argmin_{\bbmu_k \geq \bb0} L_k(\bbmu_k) :=  \argmin_{\bbmu_k \geq \bb0} \E_{\bbh_k} f(\bbmu_k, \bbh_k),  \\
&f(\bbmu_k,\bbh_k) := J (\bby_k(\bbmu_k)) + \bbmu_k^T\left( \check{\bbq}(\bbh_k,\bbp_k(\bbh_k,\bbmu_k)) - \check{\bby}_k(\bbmu_k)\right). \nonumber 
\end{align}
Typically the distribution $\ccalH_k$ is not known by the user, so the expected loss cannot be evaluated directly, but is instead replaced by an \emph{empirical risk} by taking $n$ samples labeled $\bbh_k^1, \bbh_k^2, \hdots, \bbh_k^n \in \ccalH_k^m$, (where $\bbh_k^l := [h_k^{1,l}; \hdots; h_k^{m,l}]$). In practice, such samples can be obtained through the use of short pilot signals sent from the users to measure channel conditions---see \cite{GatsisEtal14}. We then consider the empirical average loss function
\begin{equation} \label{eq_stat_av}
\hat{L}_k (\bbmu_k) := \frac{1}{n} \sum_{l=1}^n f(\bbmu_k,\bbh_k^l) := \frac{1}{n} \sum_{l=1}^n f_k^l(\bbmu_k).
\end{equation}
%
%
%
To characterize the closeness of the empirical risk $\hat{L}_k(\bbmu_k)$ with $n$ samples with respect to the expected loss $L_k(\bbmu_k)$, we define a constant $V_n$ called the \emph{statistical accuracy} of $\hat{L}_k$. The statistical accuracy $V_n$ provides a bound of the difference in the empirical and expected loss for all $\bbmu_k$ with high probability (i.e. at least $1-\gamma$ for some small $\gamma$). In other words, we define $V_n$ to be the constant that satisfies
\begin{equation}
\sup_{\bbmu_k} |\hat{L}_k(\bbmu_k) - L_k(\bbmu_k)| \leq V_n \quad \whp.
\end{equation}
The upper bounds on $V_n$ are well studied in the learning literature and in general may involve a number of parameters of the loss function $f$ as well as, perhaps most importantly, the number of samples $n$. For $\hat{L}_k(\bbmu_k)$ defined in \eqref{eq_stat_av}, a bound for the statistical accuracy $V_n$ can be obtained in the order of $\ccalO(1/\sqrt{n})$ or, in some cases, $\ccalO(1/n)$ \cite{bousquet2008tradeoffs,vapnik2013nature}. This further implies a suboptimality of $\hat{L}_k^*:=\min \hat{L}_k(\bbmu_k)$ of the same accuracy, i.e. $|L_k^* - \hat{L}_k^*| \leq 2V_n$ \cite{bousquet2008tradeoffs}. 

As is often the case in machine learning problems, the statistical accuracy informs the proper use of regularization terms in the empirical loss function. We can add regularizations to prescribe desirable properties on the empirical risk $\hat{L}_k(\bbmu_k)$, such as strong convexity, without adding additional bias beyond that already accrued by the empirical approximation. In other words, as $\hat{L}_k^*$ will be of order $V_n$ from the optimal expected value $L^*$, any additional bias of order $V_n$ or less is permissible. With that in mind, we add the regularization term $\alpha V_n/2 \|\bbmu_k\|^2$ where $\alpha>0$ to the empirical risk in \eqref{eq_stat_av} to impose strong convexity. We can further remove the non-negativity constraint on the dual variables in \eqref{eq_erm_problem} through the use of a logarithmic barrier.  To preserve smoothness for small $\bbmu_k$, we use an $\eps$-thresholded log function, defined as
\begin{align}\label{eq_log_eps}
\log_{\eps}(\bbmu_k) := \begin{cases}
\log(\bbmu_k) & \bbmu_k \geq \eps \\
\ell_{2,\eps}(\bbmu_k-\eps) & \bbmu_k < \eps,
\end{cases}
\end{align}
where $\ell_{2,\eps}(\bbmu_k)$ is a second order Taylor series expansion of $\log(\bbmu_k)$ centered at $\eps$ for some small $0<\eps<1$. We then use $-\beta V_n  \mathbf{1}^T \log_{\eps} \bbmu_k$ where $\beta>0$ as a second regularization term, and obtain a regularized empirical risk function
\begin{equation} \label{eq_stat_av_reg}
R_k (\bbmu_k) := \frac{1}{n} \sum_{l=1}^n f_k^l(\bbmu_k) + \frac{\alpha V_n}{2} \|\bbmu_k\|^2 -\beta V_n \mathbf{1}^T \log_{\eps} \bbmu_k.
\end{equation}
From here, we can seek a minimizer of the strongly convex regularized risk $R_k(\bbmu_k)$ without explicitly enforcing a non-negativity constraint on $\bbmu_k$ and find a solution with suboptimality of order $\ccalO(V_n)$ with respect to \eqref{eq_erm_problem}. Such a deterministic and strongly convex loss function as in \eqref{eq_stat_av_reg} can be minimized using a wide array of optimization methods \cite{johnson2013accelerating,defazio2014saga,mokhtari2016adaptive,eisen2017large}. However, all such methods only solve the problem for a particular epoch $k$, or otherwise assume a stationary channel distribution $\ccalH_k$ as is typical in machine learning settings.

 
\section{ERM over non-stationary channel}\label{sec_nonstat}

 The ERM problem we are interested in solving in wireless autonomous systems is further complicated by the non-stationarity of $\ccalH$, making existing solution methods insufficient. This is due to the fact that finding the minimizer to $R_k(\bbmu)$ will only provide an optimal power allocation for the respective channel distribution $\ccalH_k$. In wireless systems, we instead must \emph{continuously learn} optimal policies as the channel varies, or in other words, find optimal points for $R_k(\bbmu)$ for $k=0,1,\hdots$. To formulate the non-stationarity, however, we first define an epoch-indexed empirical risk function. While we may use a simple empirical risk as we did in \eqref{eq_stat_av}, we instead define a more general statistical loss function for a non-stationary channel using samples from the previous $M$ epochs. We define a windowed empirical loss function $\tdL_k(\bbmu)$ at epoch $k$ as
\begin{align} \label{eq_win_stat_av}
\tdL_k (\bbmu) &:= \frac{1}{M} \sum_{j=k-M+1}^k \hat{L}_j(\bbmu) 
\end{align}
By keeping a window of samples, we may retain $N=Mn$ total samples while drawing only $n$ new samples at each epoch. If the successive channel distributions $\ccalH_{k-M+1},\hdots,\ccalH_k$ are not very different, we may expect the old channel samples to still be of interest. We define the associated statistical accuracy $\tdV_N$ as the constant that satisfies
\begin{equation}\label{eq_stat_acc_w}
\sup_{\bbmu} |\tdL_k(\bbmu) - L_k(\bbmu)| \leq \tdV_N \quad \whp.
\end{equation}

Here we stress that the bounds on this constant $\tdV_N$ are not as easily obtainable or well-studied as in the stationary setting. Such a bound over non-i.i.d. samples may be dependent upon many parameters such as the sample batch size $n$, window size $M$, and correlation between successive distributions $\ccalH_{j}$ and $\ccalH_{j+1}$.  Therefore, finding precise bounds on $\tdV_N$ would require a sophisticated statistical analysis and is outside the scope of this work. We instead define a user-selected accuracy $\hat{V}$ that may estimate the statistical accuracy $\tdV_N$. We assume that $\hat{V} \geq \tdV_N$, with equality holding in cases where $\tdV_N$ is known. Using the same regularizations introduced previously, we obtain the regularized windowed empirical loss function
\begin{align} \label{eq_full_stat_av}
\tdR_k (\bbmu) &:= \frac{1}{M} \! \! \sum_{j=k-M+1}^k \!\!\!\hat{L}_j(\bbmu) + \frac{\alpha \hat{V} }{2}\|\bbmu\|^2 - \! \beta \hat{V}  \mathbf{1}^T\log_{\eps} \bbmu.
\end{align}
We subsequently define $\tbmu_k^* := \argmin_{\bbmu} \tdR_k(\bbmu)$. The definition of the loss function in \eqref{eq_full_stat_av} includes the batches of $n$ samples taken from the previous $M$ channel distributions $\ccalH_{k-M+1},\hdots,\ccalH_k$. This definition is, in a sense, a generalization of the simpler empirical risk $R_k(\bbmu)$ in \eqref{eq_stat_av_reg}. Observe that, by using a window size of $M=1$, we use only samples from the current channel and recover $R_k(\bbmu)$. In the following proposition we establish the accuracy of an optimal point of our regularized empirical risk function  $\tdR_k(\bbmu)$ relative to the optimal point of the original dual loss function $L_k(\bbmu)$.
\begin{proposition}\label{prop_order}
Consider $L_k^* = L_k(\bbmu_k^*)$ and $\tdL_k^* = \min_{\bbmu \geq \bb0} \tdL_k(\bbmu)$, and define $\tdR_k^* := \min_{\bbmu} \tdL_k(\bbmu) + \alpha \hat{V}/2 \| \bbmu\|^2 - \beta \hat{V} \mathbf{1}^T \log \bbmu$ as the optimal value of the regularized empirical risk. Define $\tdV_N$ by \eqref{eq_stat_acc_w}. Assuming $\tdV_N \leq \hat{V}$, the difference $|L_k^* - \tdR_k^*|$ is upper bounded on the order of statistical accuracy $\hat{V}$, i.e. for some $\rho>0$
\begin{equation}\label{eq_prop_diff}
| L_k^* - \tdR_k^*| \leq 2 \tdV_N + \rho \hat{V} \leq (2+\rho)\hat{V}, \quad w.h.p.
\end{equation}
\end{proposition}
\begin{myproof}
To obtain the result in \eqref{eq_prop_diff}, consider expanding and upper bounding $| L_k^* - \tdR^*| = | L_k^* - \tdL_k^*+ \tdL_k^* + \tdR_k^*| \leq  | L_k^* - \tdL_k^*| + |\tdL_k^* + \tdR_k^*|$. The first term is bounded by $2 \tdV_N$ as previously discussed. The second term, can be decomposed into the bias introduced by the logarithmic barrier $-\beta \hat{V}  \mathbf{1}^T\log \bbmu$ and the bias introduced by the quadratic regularizer $c \hat{V} /2\|\bbmu\|^2$. The former of these is known to produce an optimality bias of $(m+1) \beta \hat{V}$ \cite[Section 11.2.2]{boyd2004convex}, while the latter is known to introduce a bias on the order of $\ccalO(\hat{V})$ \cite{shalev2010learnability}. Combining these, we get a total suboptimality between the regularized risk function optimal and the true optimal of $2 \tdV_N + \rho \hat{V}$ for some constant $\rho>0$. As we assume that $\tdV_N \leq \hat{V}$, the rightmost bound in \eqref{eq_prop_diff} follows.
\end{myproof}


A key observation to be made here is that any exact solution to \eqref{eq_full_stat_av} only minimizes the expected loss $L_k$ to within accuracy $\hat{V}$ (assuming $\hat{V} \geq \tdV_N$). There is therefore no need to minimize \eqref{eq_full_stat_av} exactly but is in fact sufficient to find a $\hat{V}$-accurate solution, as this incurs no additional error relative to the statistical approximation itself. While many optimization methods can be used to find a minimizer to \eqref{eq_full_stat_av}, we demonstrate in the next section that fast second order methods can be used to learn \emph{approximate} minimizers---and by Proposition \ref{prop_order} approximately solve \eqref{eq_erm_problem}---at each epoch $k$ with just single updates as the channel distribution $\ccalH_k$ changes, thus tracking near-optimal points at every epoch. This is done by exploiting an important property of second order optimization methods, namely \emph{local quadratic convergence}. 

\begin{remark}\label{remark_log}\normalfont
Observe in the text of Proposition \ref{prop_order} that we define $\tdR_k^*$ to be the optimal point of the loss function $\tdL_k(\bbmu_k)$ regularized with a standard log barrier $-\log(\bbmu_k)$, rather than the thresholded barrier $-\log_{\eps}(\bbmu_k)$ used in the definition in \eqref{eq_full_stat_av}. Indeed, using the thresholded barrier does not explicitly enforce nonnegativity for values smaller than $\eps$. However, this thresholding is necessary to preserve smoothness of the barrier, which will be necessary for the proof of Lemma \ref{main_lemma} in Section \ref{sec_convergence}. The threshold $\eps$ can be made as small as necessary to enforce nonnegativity, although this comes at the cost of a worse smoothness constant. In practice, however, we observe this thresholding to not be explicitly needed and is just included here for ease of analysis. We also stress that the smoothness constant itself does not play a pivotal role in the proceeding analysis.
\end{remark}

\subsection{Learning via Newton's Method}\label{sec_ada_newton}
In this paper, we use Newton's method to approximately minimize \eqref{eq_full_stat_av} efficiently as the channel $\ccalH_k$ changes over epochs. Motivated by the recent use of Newton's method in solving large scale ERM problems through adaptive sampling policies \cite{mokhtari2016adaptive,eisen2017large}, we use the $N$ samples drawn from recent distributions to find an iterate $\bbmu_k$ that approximately solves for $\tbmu_k^*$. At the next epoch, the iterate $\bbmu_{k}$ provides a ``soft'' start towards finding a point $\bbmu_{k+1}$ that approximately minimizes $\tdR_{k+1}(\bbmu)$. In this way, with single iterations we may find near-optimal solutions for each regularized empirical loss function, and thereby efficiently learn the optimal power allocation of the wireless channel as the channel distribution evolves over time epochs.

We proceed by presenting the details of Newton's method. At epoch $k$, we compute a new iterate $\bbmu_{k+1}$ by subtracting from the current iterate $\bbmu_k$ the product of the Hessian inverse and the gradient of the function $\tdR_{k+1}(\bbmu_k)$. For the empirical dual loss function $\tdR_{k}$ defined in \eqref{eq_full_stat_av}, we define the gradient $\nabla \tdR_k(\bbmu)$ and Hessian $\nabla^2 \tdR_{k}(\bbmu)$. 
The new approximate solution $\bbmu_{k+1}$ is then found from current approximate solution $\bbmu_k$ using the Newton update
\begin{equation}\label{eq_update}
\bbmu_{k+1} = \bbmu_{k} - \bbH_{k+1}^{-1} \nabla \tdR_{k+1}(\bbmu_{k}),
\end{equation}
where we use $\bbH_{k+1} := \nabla^2 \tdR_{k+1}(\bbmu_k)$ as simplified notation.

To understand the full algorithm, consider that $\bbmu_{k}$ is a $\hat{V}$-accurate solution of current loss function $\tdR_{k}$, i.e. $\tdR_{k}(\bbmu_{k}) - \tdR_{k}^* \leq \hat{V}$. Recall that the new loss function $\tdR_{k+1}$ differs from $\tdR_k$ only in the discarding of old samples $\hat{L}_{k-M+1}$ and inclusion of samples $\hat{L}_{k+1}$ drawn from $\ccalH_{k+1}$. If we consider that the distributions are varying slowly across successive time epochs, i.e. $\ccalH_{k+1}$ is close to $\ccalH_{k}$, then the respective loss functions $\tdR_{k+1}$ and $\tdR_k$ and their optimal values $\tdR_{k+1}^*$ and $\tdR_k^*$ will also not differ greatly under some smoothness assumptions. Therefore, under such conditions a single step of Newton's method as performed in \eqref{eq_update} can in fact be sufficient to reach a $\hat{V}$-accurate solution of the new loss function $\tdR_{k+1}$. This is possible precisely because of the Newton method's property of \emph{local quadratic convergence}, meaning that Newton's method will find a near-optimal solution very quickly when it is already in a local neighborhood of the optimal point. Given then a $\hat{V}$-accurate solution $\bbmu_0$ of initial loss $\tdR_{0}$, the proceeding and all subsequent iterates $\bbmu_k$ will remain within the statistical accuracy of their respective losses $\tdR_k$ as the channel distribution varies over time. The formal presentation of the exploitation of this property and other technical details of this result are discussed in Section \ref{sec_convergence} of this paper.

The learning algorithm is presented in Algorithm \ref{alg:AdaNewton}. After preliminaries and initializations in Steps 1-4, the backtracking loop starts in Step 5. Each iteration begins in Step 6 with the the drawing of $n$ samples from the new channel distribution $\ccalH_{k}$ and discarding of old samples from $\ccalH_{k-M}$ to form $\tdR_{k}$. Note that samples will be only be discarded for $k>M$. The gradient $\nabla \tdR_{k}$ and Hessian $\bbH_{k}$ of the regularized dual loss function are computed in Step 7. The Newton step is taken with respect to $\tdR_{k+1}$ in Step 8. In Step 9, the optimal primal variables are computed with respect to the updated dual variables. This includes both the auxiliary variables $\bby(\bbmu_{k})$ and the power allocation policy $\bbp(\bbh,\bbmu_{k})$ itself. Because there are function and channel system parameters that are not known in practice, we include a backtracking step for the parameters $n$ and $M$ in Step 10 to ensure the new iterate $\bbmu_{k}$ is within the intended accuracy $\hat{V}$ of $\bbmu_{k}^*$. Further details on the specifics of the backtracking procedure are discussed in Section \ref{sec_implementation} after the presentation of the theoretical results.

%
{\linespread{1.3}
\begin{algorithm}[t] \begin{algorithmic}[1]\small
\STATE \textbf{Parameters:} Sample size increase constants $n_0>0$, $M_0 \geq 1$  backtracking params $0<\gamma<1<\Gamma$, and accuracy $\hat{V}$.
\STATE \textbf{Input:} Initial sample size $n=m_0$ and 
                       argument $\bbmu_{n} = \bbmu_{m_0}$ with 
                       $\| \nabla \tdR_{n}(\bbmu_{k+1})\| < (\sqrt{2 \alpha}) \hat{V}$
\FOR [main loop]{$k = 0,1,2,\hdots$} 
   \STATE Reset factor $n=n_0$, $M = M_0$ .     
   \REPEAT  [sample size backtracking loop] 
      \STATE Draw $n$ samples from $\ccalH_{k}$, discard from $\ccalH_{k-M}$. 
      \STATE Compute Gradient $\nabla \tdR_{k} (\bbmu_{k-1})$, Hessian $\bbH_{k}$. \ 
      \STATE Newton Update [cf. \eqref{eq_update}]: 
        $$
                 \bbmu_{k} = \bbmu_{k-1} - \bbH_{k}^{-1} \nabla \tdR_{k} (\bbmu_{k-1})
             $$
      \STATE Determine power allocation, aux. variables [cf. \eqref{eq_primal_recovery}, \eqref{eq_primal_recovery_2}]:  
       \begin{align}
       		p_k^i(h_k^i,\bbmu_k) &= \argmax_{p_k^i \in [0,p_0]}  \mu_k^i q(h_k^i,p_k^i(h_k^i)) -  \tdmu p_k^i(h_k^i) ,  \nonumber \\
		y_k^i(\bbmu_k) &= \argmax_{y_k^i}  J^i (y_k^i) - \mu_k^i y_k^i. \nonumber
	\end{align}
      \STATE Backtrack sample draw $n=\Gamma n$, window size $M = \gamma M$.       
   \UNTIL {$\| \nabla \tdR_{k}(\bbmu_{k})\| < (\sqrt{2 \alpha}) \hat{V}$} 
\ENDFOR
\end{algorithmic}

\caption{Learning via Newton's Method}\label{alg:AdaNewton} \end{algorithm}}

\section{Convergence Analysis}\label{sec_convergence}
 
 In this section we provide a theoretical analysis of the Newton learning update in \eqref{eq_update}. We do so by first analyzing the convergence properties of the ERM problem in \eqref{eq_full_stat_av}. We subsequently return to the WCP in \eqref{eq_power_problem_y} and establish a control performance result.

 \subsection{Convergence of ERM problem}\label{sec_erm_conv}
 We begin by analyzing the ERM formulation of the power allocation problem in \eqref{eq_full_stat_av} and establish a theoretical result that, under certain conditions, guarantees each iterate $\bbmu_{k}$ is within the statistically accuracy of the risk function at epoch $k$. Our primary theoretical result characterizes such conditions dependent on statistical accuracy and rate of non-stationarity. We begin by presenting a series of assumptions made in our analysis regarding the dual loss functions $f$.
\begin{assumption}\label{ass_convexity}
The expected loss function $L_k$ and empirical loss functions $f(\bbmu, \bbh_k)$ are convex with respect to $\bbmu$ for all values of $\bbh_k$. Moreover, their gradients $\nabla L_k(\bbmu)$ and $\nabla f(\bbmu,\bbz)$ are Lipschitz continuous with constant $\Delta$. 
\end{assumption}


\begin{assumption}\label{ass_self_concor}
The loss functions $f(\bbmu,\bbh)$ are self-concordant with respect to $\bbmu$ for all $\bbh$, i.e. for all $i$, $$| \partial^3/\partial \mu_i^3 f(\bbmu,\bbh)| \leq 2 \partial^2/\partial \mu_i^2 f(\bbmu,\bbh)^{3/2}.$$
\end{assumption}

Assumption \ref{ass_convexity} implies that the regularized empirical risk gradients $\nabla \tdR_{k}$ are Lipschitz continuous with constant $\Delta+c\hat{V}$ where $c:=\alpha + \beta/\eps^2$ and $\alpha, \beta, \eps$ are the regularization constants in \eqref{eq_full_stat_av}. The function $\tdR_k$ is also strongly convex with constant $\alpha\hat{V}$. This implies an upper and lower bound of the eigenvalues of the Hessian of $\tdR_k$, namely
\begin{equation}\label{eq_hessian_bounds}
\alpha \hat{V} \bbI \preceq \bbH_{k} \preceq (\Delta + c\hat{V})\bbI.
\end{equation}
Assumption \ref{ass_self_concor} states the loss functions are additionally self concordant, which is a common assumption made in the analysis of second-order methods---see, e.g. \cite[Ch. 9]{boyd2004convex}, for such an analysis. It also follows that the functions $\tdR_{k+1}$ are therefore self concordant because both the quadratic and thresholded log regularizers are self-concordant.  We present a brief remark regarding the implications of these assumptions on the dual risk function  on the wireless control problem.

\begin{remark}\label{remark_assumption}\normalfont
We state the preceding assumptions in terms of the sampled dual functions $f$ due to their direct use in the proceeding analysis. However, they indeed have implications on the primal domain problem in \eqref{eq_power_problem_y}. While the dual function is always convex, the smoothness condition in Assumption \ref{ass_convexity} can be obtained from the strong concavity of the control performance $\sum_i J^i$ with strong concavity $1/\Delta$.  The self-concordance property on the dual function in Assumption \ref{ass_self_concor}, however, is not easily derived from properties of $J^i(\cdot)$ or $q(\cdot)$. We point to work that establishes self concordance of the dual for various machine learning problems \cite{owen2013self,necoara2009interior}.  
\end{remark}
%
 The two preceding assumptions deal specifically with the properties of the empirical dual loss functions used in the ERM problem. To connect the solving of the sampled functions $f^l$ with the expected loss function $L$, we additionally include two assumptions regarding the statistics of the expected and empirical losses.
\begin{assumption}\label{ass_grad_cond}
The difference between the gradients of the empirical risk $\hat{L}_{k}$ and the statistical average loss $L_k$ is bounded by $V_{N}^{1/2}$ for all $\bbmu$ and $k$ with high probability,
\begin{align}\label{eqn_loss_minus_erm_2}
   \sup_{\bbmu}\|\nabla L_k(\bbmu) - \nabla \hat{L}_{k}(\bbmu) \|  \leq V_{N}^{1/2},  \qquad\whp.
\end{align}
\end{assumption}
\begin{assumption}\label{ass_nonstat}
The difference between two successive expected loss functions $L_k(\bbmu) = \E_{h_k} f(\bbmu, h_k)$ and $L_{k+1}(\bbmu) = \E_{h_{k+1}} f(\bbmu, h_{k+1})$ and the difference between their gradients are bounded respectively by a bounded sequence of constants $\{D_k\},\{\bar{D}_k\}\geq 0$ for all $\bbmu$,
\begin{align}\label{eqn_nonstat_bound}
   \sup_{\bbmu} | L_k(\bbmu) - L_{k+1}(\bbmu) |  \leq D_k,
\end{align}
\begin{align}\label{eqn_nonstat_bound_grad}
   \sup_{\bbmu}\| \nabla L_k(\bbmu) - \nabla L_{k+1}(\bbmu) \|  \leq \bar{D}_k.
\end{align}
\end{assumption}
Assumption \ref{ass_grad_cond} bounds the difference between gradients of the expected loss and the empirical risk with $N$ samples by $V_{N}^{1/2}$, which can be readily obtained using the law of large numbers. Assumption \ref{ass_nonstat} bounds the point-wise difference in the expected loss functions and their gradients at epochs $k$ and $k+1$. This can be interpreted as the rate at which the channel evolves between epochs, and is used to establish that optimal dual variables for two consecutive empirical risk functions $\tdR_{k}$ and $\tdR_{k+1}$ are not very different. We discuss practical implications of this assumption in Section \ref{sec_implementation}.

\begin{remark}\label{remark_ass4}\normalfont
Observe that the bounds provided in Assumption \eqref{ass_nonstat} are with respect to the dual function rather than explicitly on the non-stationary statistics of the channel. They are provided as such because this is the manner in which the non-stationarity appears in the proceeding analysis. To see how the channel characteristics play a role in the provided bound, consider that, e.g., \eqref{eqn_nonstat_bound} can be expanded using the definition of the dual function $L_k(\bbmu)$ as
\begin{align}\label{eq_neww}
&\sup_{\bbmu} 
|\E_{\bbh_k} \left\{\max_{\bbp \in \ccalP}  \bbmu^T\check{\bbq}(\bbh_k,\bbp(\bbh_k)) \right\} \\ &\qquad - \E_{\bbh_{k+1}} \left\{\max_{\bbp \in \ccalP}  \bbmu^T\check{\bbq}(\bbh_{k+1},\bbp(\bbh_{k+1})) \right\}|   \leq D_k. \nonumber
\end{align}
The exact condition this imposes upon the channel distribution variation thus depends both on the form of the distributions $\ccalH_k$, $\ccalH_{k+1}$, and the function $q(\bbh,\bbp)$. Thus, the exact manner in which the varying channel conditions effect this bound are indeed problem-specific, and a generic condition on non-stationarity of the channel is only present in the proceeding analysis indirectly through the condition in \eqref{eq_neww}.
\end{remark}

The proceeding analysis is organized in the following manner. Our goal is to establish conditions on the parameters of the statistical accuracy---$\hat{V}$---and the non-stationarity---$D_k$ and $\bar{D}_k$---that guarantee that, starting from an approximate solution to $\tdR_k$, a single step of Newton's method generates an approximately accurate solution to $\tdR_{k+1}$. From there, we can recursively say that, assuming an initial point $\bbmu_0$ that is within the intended accuracy of $\tdR_0$, the method will continue to find a $\hat{V}$-accurate solution at each epoch as the channel distribution changes with $k$. We achieve this result in two steps. We first find a condition that guarantees that a $\hat{V}$-accurate solution of $\tdR_k$ is also in the quadratic convergence region of $\tdR_{k+1}$. Second, we find a condition that guarantees that such a point within the quadratic convergence region of $\tdR_{k+1}$ will reach its intended accuracy with a single update as in \eqref{eq_update}. 

We begin by establishing the condition in the first step, namely that a $\hat{V}$-accurate solution to $\tdR_k$, labeled $\bbmu_k$ is in in the quadratic convergence region of $\tdR_{k+1}$ if certain conditions hold. The quadratic convergence region is a region local to the optimum in which Newton's method is known to converge at a fast quadratic rate. The analysis of Newton's method commonly characterizes quadratic convergence in terms of a quantity called the Newton decrement, explicitly defined as $\lambda_{k+1} (\bbmu) := \|\nabla^2 \tdR_{k+1}(\bbmu)^{-1/2}  \nabla \tdR_{k+1}(\bbmu) \|$. We say the dual iterate $\bbmu$ is in the quadratic convergence region of $\tdR_{k+1}$ when $\lambda_{k+1}(\bbmu) < 1/4$---see \cite[Chapter 9.6.4]{boyd2004convex}. In the following proposition, we give conditions under which any iterate $\bbmu_k$ that is within the accuracy $\hat{V}$ of the optimal point $\tdR_k^*=\min_\mu \tdR_k(\mu)$ is also within the quadratic convergence region of the subsequent loss function $\tdR_{k+1}$. 

\begin{lemma}\label{main_lemma}
Consider $\bbmu_{k}$ as a $\hat{V}$-accurate optimal solution of the loss $\tdR_{k}$, i.e., $\tdR_{k}(\bbmu_k)-\tdR_{k}^* \leq \hat{V}$. In addition, define $\lambda_{k+1} (\bbmu):=\left(\nabla \tdR_{k+1}(\bbmu) ^T \nabla^2 \tdR_{k+1}(\bbmu) ^{-1} \nabla \tdR_{k+1}(\bbmu) \right)^{1/2}$ as the Newton decrement of variable $\bbmu$ associated with the loss $\tdR_{k+1}$. If Assumptions \ref{ass_convexity}-\ref{ass_nonstat} hold, then Newton's method at point $\bbmu_k$ is in the quadratic convergence phase for the objective function $\tdR_{k+1}$, i.e., $\lambda_{k+1}(\bbmu_k)<1/4$, if we have
\begin{equation}\label{prop_result}
 \left( \frac{2(\Delta+c\hat{V})\hat{V}}{\alpha\hat{V}}\right)^{1/2} + \frac{2 \tdV_N^{1/2} + \bar{D}_k}{(\alpha \hat{V})^{1/2}} < \frac{1}{4}.
\quad \text{w.h.p.}
\end{equation}
\end{lemma}
\begin{myproof}
See Appendix.
\end{myproof}

Lemma \ref{main_lemma} provides the first necessary condition in our analysis by identifying the statistical parameters under which every iterate $\bbmu_k$ is in the quadratic region of $\tdR_{k+1}$. From here we can show the second step, in which such a point in the quadratic convergence region of $\tdR_{k+1}$ can reach its statistical accuracy with a single Newton step as given in \eqref{eq_update}. To achieve this, we first present the following lemma that upper bounds the sub-optimality of the point $\bbmu _k$ with respect to the optimal solution of $R^*_{k+1}$.
 
%
\begin{lemma}\label{lemma_sub_m}
Consider a point $\bbmu_k$ that minimizes the loss function $\tdR_k$ to within accuracy $\hat{V}$, i.e. $\tdR_{k}(\bbmu_k) - \tdR_{k}^* \leq \hat{V}$. Provided that Assumptions \ref{ass_convexity}-\ref{ass_nonstat} hold, the sub-optimality $\tdR_{k+1}(\bbmu_k) - \tdR_{k+1}^*$ is upper bounded w.h.p. as
\begin{align} \label{eq_lemma_sub_m_2}
\tdR_{k+1}(\bbmu_k) - \tdR_{k+1}^* \leq 4\tdV_N + \hat{V}+ 2D_k
\end{align}
\end{lemma}
%
\begin{myproof}
See Appendix.
\end{myproof}

In Lemma \ref{lemma_sub_m} we establish a bound on the suboptimality of $\bbmu_k$ with respect to $\tdR_{k+1}$. The following lemma now bounds the suboptimality of $\bbmu_{k+1}$ in terms of the suboptimality of $\bbmu_k$ with a quadratic rate.

\begin{lemma}\label{quadratic_convg_lemma}
Consider $\bbmu_{k}$ to be in the quadratic neighborhood of the loss $\tdR_{k+1}$, i.e., $\lambda_{k+1} (\bbmu_k)\leq1/4$. Recall the definition of the variable $\bbmu_{k+1}$ in \eqref{eq_update} as the updated variable using Newton's method. If Assumptions \ref{ass_convexity}-\ref{ass_grad_cond} hold, then the difference $\tdR_{k+1}(\bbmu_{k+1})-\tdR_{k+1}^*$ is upper bounded by
\begin{equation}\label{prop_result_2}
		\tdR_{k+1}(\bbmu_{k+1})-\tdR_{k+1}^* \leq 144 (\tdR_{k+1}(\bbmu_{k})-\tdR_{k+1}^* )^2.
\end{equation}
\end{lemma}
\begin{myproof}
See Appendix.
\end{myproof}

With Lemma \ref{quadratic_convg_lemma} we establish the known quadratic rate of convergence of the suboptimality of the Newton update in \eqref{eq_update}. Now by substituting the upper bound on $\tdR_{k+1}(\bbmu_{k})-\tdR_{k+1}^*$ from Lemma \ref{lemma_sub_m}, a condition can easily be derived under which the suboptimality of the new iterate is within the accuracy $\hat{V}$ of $\tdR_{k+1}$. Using the results of Lemmata \ref{main_lemma}-\ref{quadratic_convg_lemma}, we present our main result in the following theorem.

%
\begin{theorem}\label{theorem_main_result}
Consider Newton's method defined in \eqref{eq_update} and the full learning method detailed in Algorithm \ref{alg:AdaNewton}. Define $\tdV_N$ to be the statistical accuracy of $N=Mn$ samples  by \eqref{eq_stat_acc_w}, with $n$ samples taken from each of the $M$ most recent channel distributions $\ccalH_k$. Further consider the variable $\bbmu_{k}$ as a $\hat{V}$-optimal solution of the loss $\tdR_{k}$, and suppose Assumptions \ref{ass_convexity}-\ref{ass_nonstat} hold. If the sample size $n$ and window size $M$ are chosen such that the following conditions
\begin{equation}\label{cond_1}
 \left( \frac{2(\Delta+c\hat{V})\hat{V}}{\alpha\hat{V}}\right)^{1/2} + \frac{2\hat{V}^{1/2} + \bar{D}_k}{(\alpha \hat{V})^{1/2}} < \frac{1}{4}
\end{equation}
\begin{equation}\label{cond_2}
144 (4\tdV_N + \hat{V} + 2D_k)^2
\leq \hat{V}
\end{equation}
are satisfied, then the variable $\bbmu_{k+1}$ computed from \eqref{eq_update} has the suboptimality of $\hat{V}$ with high probability, i.e., 
\begin{equation}\label{imp_result}
 \tdR_{k+1}(\bbmu_{k+1})- \tdR_{k+1}^* \leq \hat{V}, \qquad \whp.
\end{equation}
 \end{theorem}
 
The inequalities \eqref{cond_1}-\eqref{cond_2} in Theorem \ref{theorem_main_result} specify conditions under which $\bbmu_{k+1}$ as generated by \eqref{eq_update} is a $\hat{V}$-optimal solution of $\tdR_{k+1}$. Note that these conditions come directly from the preceding lemmata. Thus, when these conditions are satisfied, single iterations of Newton's method at each epoch $k$---as detailed in Algorithm \ref{alg:AdaNewton}---successively generate approximately optimal dual parameters. A further discussion of the satisfaction of such conditions in regards to practical implementation is provided later in Section \ref{sec_implementation}. We first extend the theoretical result of Theorem \ref{theorem_main_result} to establish properties of the resulting WCP solution.

\begin{remark}\label{remark_theorem}\normalfont
Observe in Theorem \ref{theorem_main_result} that the provided conditions cannot be satisfied if the true statistical accuracy $\tdV_N$ is greater than the selected $\hat{V}$. While we assume in our analysis this is not the case, (i.e. $\hat{V}$ is a conservative estimate of $\tdV_N$), this may not be guaranteed if very little information is known about $V_N$. In the case $\hat{V} < V_N$, we point out that the results in Theorem 1 can simply be modified by replacing achieved accuracy $\hat{V}$ by $V_N$. In other words, the accuracy we can achieve is limited by the greater of these terms. We do not go through the details of this analysis for clarity of presentation, but such result can be obtained through the same steps of the preceding analysis. 
\end{remark}


\subsection{Sub-optimality in wireless control system}\label{sec_subopt}

Because the proposed Newton method indeed solves \eqref{eq_erm_problem} to within a statistical approximation $\hat{V}$, it is important to consider the effect of such an approximation on the original WCP in \eqref{eq_power_problem_y}. In this section we provide a sequence of results that characterize the accuracy of the solutions generated by the Newton update in \eqref{eq_update} in the original primal control problem in \eqref{eq_power_problem_y}. Firstly, recall the constraints in \eqref{eq_power_problem_y} reflect both a power budget limited by $p_{\max}$ and that the auxiliary variable $y^i$ should not exceed the expected packet success function $q(\cdot)$. In solving the dual problem approximately, we may then also violate these constraints by a small margin.  We can specifically characterize such a constraint violation, as well as address the suboptimality in terms of the primal objective. Both these results together can then be combined to demonstrate the stability of the switched system WCP introduced in Example \ref{example1}. To do so, we first introduce an assumption regarding the feasibility and boundedness of the dual loss solutions $L_k^*$ and the optimal dual point $\bbmu_k^*$.

\begin{assumption}\label{ass_feasible}
For all epochs $k$, the problem in \eqref{eq_power_problem_y} under distribution $\ccalH_k$ is strictly feasible. There also exists constants $\ccalK$ and $\hat{\ccalK}$ such that the optimal dual objective value $L_k^*$ is bounded as $L_k^* \leq \ccalK$ and optimal dual variable bounded as $\|\bbmu_k^*\| \leq \hat{\ccalK}$.
\end{assumption}

From strict feasibility of the primal problem in \eqref{eq_power_problem_y}, we also obtain a finite upper bound on the value of the dual function. This can be used with the suboptimality result in Theorem \ref{theorem_main_result} to bound the norm of the dual variables $\bbmu_k$ generated from the Newton update in \eqref{eq_update}. This is presented in the following corollary.
\begin{corollary}\label{corollary_primal_bound}
The norm of the dual variables $\bbmu_k$ generated by the update in \eqref{eq_update} is bounded as $\|\bbmu_k\| \leq \sqrt{(2/\alpha)} + \hat{\ccalK}$.
\end{corollary}
\begin{myproof}
From strong convexity we have that $\|\bbmu_k - \tbmu_k^*\|^2 \leq (2/\alpha \hat{V}) (\tdR_k(\bbmu_k) - \tdR_k^*)$. Using the reverse triangle inequality with \eqref{imp_result} and Assumption \ref{ass_feasible}, we obtain the intended result. 
\end{myproof}

Observe that the boundedness of the solutions to the regularized dual function in Assumption \ref{ass_feasible} in effect states that, for all distributions $\ccalH_k$, the empirical, or sampled, versions of the constrained problem in \ref{eq_power_problem_y} will be strictly feasible. From here, we can establish through duality a bound on each constraint violation that may occur from solving the dual problem to its statistical accuracy. This result is stated in the following proposition.

\begin{proposition}\label{feasible_error_prop}
Consider $\bbmu_{k}$ to be a $\hat{V}$-optimal minimizer of $\tdR_k$, i.e. $\tdR_k(\bbmu_k) - \tdR_k^* \leq \hat{V}$. Further consider $\bbp(\bbh,\bbmu_k)$ and $\bby(\bbmu_k)$ to be the Lagrangian maximizers over dual parameter $\bbmu_k$. If Assumptions \ref{ass_convexity} and \ref{ass_feasible} hold, then the norm of the constraint violations in \eqref{eq_power_problem_y} can each be upper bounded as
\begin{align}\label{eq_j1}
		\left| \sum_{i=1}^m \E_{h_k^i} (p^i(h_k^i,\bbmu)) - p_{\max} \right| &\leq \sqrt{2\Delta(\tdV_N + C \hat{V})}, \\
		\left\|\bby(\bbmu_k) -  \E_{\bbh_k} \left\{\bbq(\bbh_k,\bbp(\bbh_k,\bbmu_k))\right\} \right\| &\leq \sqrt{2\Delta(\tdV_N + C \hat{V})}, \label{eq_j2}
\end{align}
where $C := 1 + \rho + \beta \kappa$ and $\kappa$ such that $\mathbf{1}^T \log_{\eps} (\bbmu_k) \leq \kappa$.
\end{proposition}

\begin{myproof}
See Appendix.
\end{myproof}

In Proposition \ref{feasible_error_prop}, we establish a bound that is proportional to $\hat{V}$ on the violation of the constraints in \eqref{eq_power_problem_y}. There are two points to be stressed here. First, is that this constraint violation can indeed be made small by controlling the target accuracy $\hat{V}$. Additionally, we point out that the violation of the budget constraint can be controlled by adding a slack term to the maximum power as $\hat{p}_{\max} = p_{\max} - 2\Delta C \hat{V}$. In this way, any such violation will still be within the true intended budget $p_{\max}$.

We proceed by establishing suboptimality of the generated variables $\bby(\bbmu_k)$ in terms of control performance. Recall the final result in Theorem \ref{theorem_main_result} that establishes at each epoch $k$, the current dual function value $\tdR_k(\bbmu_k)$ will be within accuracy $\hat{V}$ of the optimal value $\tdR_k(\tbmu_k^*)$ (after satisfying the necessary conditions). To establish that the control systems induced by such dual parameters $\bbmu_k$ remain stable, we first connect the accuracy of the dual function value to the accuracy of associated primal variables $\bbp(\bbh,\bbmu_{k})$ and $\bby(\bbmu_{k})$ with respect to their optimal values $\bbp_k^*(\bbh) := \bbp(\bbh,\bbmu^*_{k})$ and $\bby_k^* := \bby(\tbmu_k^*)$. This bound is established in the following theorem.

\begin{theorem}\label{cor_primal_error}
Consider $\bbmu_{k}$ to be a $\hat{V}$-optimal minimizer of $\tdR_k$, i.e. $\tdR_k(\bbmu_k) - \tdR_k^* \leq \hat{V}$. Further consider $\bbp(\bbh,\bbmu_k)$ and $\bby(\bbmu_k)$ to be the Lagrangian maximizers over dual parameter $\bbmu_k$. Under Assumptions \ref{ass_convexity}-\ref{ass_feasible} the primal objective function sub-optimality $J(\bby(\bbmu_k)) - J(\bby_k^*)$ can be upper bounded as
\begin{equation}\label{eq_lemma_primal_error}
		J(\bby(\bbmu_k)) - J(\bby_k^*) \leq (1+C)\Delta\left(\frac{1}{\alpha} + 2 \hat{V} (\sqrt{2/\alpha} + \hat{\ccalK})\right).
\end{equation}
\end{theorem}

\begin{myproof}
See Appendix.
\end{myproof}

In Theorem \ref{cor_primal_error}, we derive a bound on the suboptimality of the primal objective function $J(\bby)$ that is proportional also to the statistical accuracy $\hat{V}$ plus a constant. Recall that this function is, in general, a measure of the control performance of the system. Thus, solving the dual problem approximately indeed can be translated into approximate accuracy in terms of our original utility metric with respect to the control system. In many problems, the performance $J(\bby)$ will also effectively establish a stability margin for control systems that have unstable regions of operation. To demonstrate the effect of using the proposed Newton's method over a non-stationary wireless channel, we return to the switched dynamical system in Example \ref{example1}.

\subsection{Stability of switched dynamical system (Example \ref{example1})}\label{sec_stability}
Consider the switched dynamical system given in \eqref{eq:system} and the derived performance metric $J(\bby)$ in \eqref{eq_j_ex1} that tracks the asymptotic behavior of the state $x_t$. In this system, if the open loop gain is unstable $|A_o|>1$ it can indeed cause the system to grow in an unstable manner if the system is not closed sufficiently often. As mentioned in Example \ref{example1} the system reaches instability if $y A_c^2 + (1-y) A_o^2$ becomes close to $1$. A question of interest in this example is, using the power allocation policy found using Newton's method over a time-varying channel, whether or not the system remains stable over time. We can indeed demonstrate this to be true with the following argument.

From Theorem \ref{cor_primal_error}, we obtained that the primal suboptimality with respect to the control performance function $J(\bby)$ is bounded by a term proportional to $\hat{V}$. Assuming that $J(\bby_k^*)$ is finite for all epochs $k$, it follows then that the generated performance $J(\bby(\bbmu_k))$ is also finite. Considering the expression for $J^i(y^i)$ given in \eqref{eq_j_ex1}, this is finite if and only if the denominator is positive, i.e., there exists a $\omega$ such that
\begin{equation}\label{eq_y_bound}
1 - y^i(\bbmu_k)( (A^i_o)^2- (A^i_c)^2)  \leq \omega < 1
\end{equation}
 at all epochs $k$.

Moreover from Proposition~\ref{feasible_error_prop} we also have that the actual packet success rate during epoch $k$ satisfies
\begin{align}
\E_{h_k} \left\{q(h_k^i,p(h_k^i,\mu_k))\right\} \geq y^i(\bbmu_k) - \sqrt{2\Delta(\tdV_N + C \hat{V})},
\end{align}
If the statistical accuracy at the right hand side of this expression is sufficiently small, then using \eqref{eq_y_bound} we also get that
\begin{equation}\label{eq_packet_rate_bound}
1 - \E_{h_k} \left\{q(h_k^i,p(h_k^i,\mu_k))\right\} ((A^i_o)^2 - (A^i_c)^2) \leq \tilde\omega < 1
\end{equation}
In particular this holds if $\sqrt{2\Delta(\tdV_N + C \hat{V})} ({(A^i_o)^2 - (A^i_c)^2})< {1 - \omega}$.

Substituting \eqref{eq_packet_rate_bound} back into the recursive expression in \eqref{eq_system_d}, we get that the variance of the state at each time step satisfies
\begin{equation}\label{eq_system_d2}
	\mathbb{E}(x^i_{t+1})^2 \leq \tilde\omega\mathbb{E}(x^i_t)^2 + W^i.
\end{equation} 
 Operating recursively and using the geometric series as in Example~\ref{example1}, we can bound \eqref{eq_system_d2} as
\begin{align}
	\mathbb{E}(x^i_{t+1})^2 &\leq 
	 \tilde\omega^{t+1} \mathbb{E}(x^i_0)^2 + W^i\frac{1-\tilde\omega^{t+1}}{1-\tilde\omega} .\label{eq_system_d3}
\end{align} 
As both terms on the right hand side of \eqref{eq_system_d3} are finite, we can conclude that the state variables remain bounded in variance for all $t$ in the non-stationary channel.

\section{Details of Implementation}\label{sec_implementation}

In this section we provide a discussion of necessary considerations for practical implementation of the result in Theorem \ref{theorem_main_result}. Observe that the conditions in \eqref{cond_1} and \eqref{cond_2} are functions of four primary terms, $\hat{V}$, $\tdV_N$, $D_k$, and $\bar{D}_k$. While $\hat{V}$ is user-selected, the latter three terms come directly from statistical properties of the control performance functions and the channel distribution. They can, however, be indirectly controlled for with some careful implementation techniques.

First, consider that the latter two terms $D_k$  and $\bar{D}_k$ provide a bound on the difference of the neighboring expected loss functions $L_k$ and $L_{k+1}$ and their gradients, respectively. Thus, these terms collectively can be interpreted as a bound on the degree of non-stationarity of the channel distribution $\ccalH$ between successive time iterations, or in other words the rate at which the channel changes over time epochs. In a practical sense, this rate is controllable by determining how much real time makes up a single discrete time epoch. That is, time epochs $k$ and $k+1$ that are closer together in a real time-sense will naturally have a lower bound for $D_k$, and $\bar{D}_k$, assuming the rate of change of the channel distribution is indeed smooth. In this sense, $D_k$  and $\bar{D}_k$ can be lowered to satisfy the conditions in \eqref{cond_1} and \eqref{cond_2} by considering shorter time between discrete epochs. This is to say that, because the channel conditions are not in our control, if necessary we may change the rate at which we apply our algorithm in a real time sense. By using shorter epochs, we collect channel samples and run the proposed Newton step more often to adapt to quickly changing channel conditions.
 
The second term present in the conditions of Theorem \ref{theorem_main_result}---namely $\tdV_N$---represents the statistical accuracy of the non-i.i.d. samples taken from the window of $M$ most recent channel distributions with respect to the current channel distribution. A condition on $\tdV_N$ in fact then indirectly provides conditions on the sample size $n$ and window size $M$ used to define $\tdR_{k}$ necessary to learn a $\hat{V}$-optimal solution. We reiterate here that, in the simpler setting of $M=1$, a well-studied bound on $\tdV_N$ exists of the order $\ccalO(1/\sqrt{n})$. For the case of windowed sampling the bound on $\tdV_N$ can nonetheless still be varied through various choices of window size $M$ and sample draw size $n$. However, because the exact nature of both $\tdV_N$ and $D_k$ come from statistical properties not known in practice, precise selection of such parameters $n$ and $M$ can be chosen via a standard backtracking procedure. 

The details of the backtracking procedure can be seen in Steps 10 and 11 in Algorithm \ref{alg:AdaNewton}. At each epoch $k$, the parameters $n$ and $M$ are initialized to $n_0$ and $M_0$ in Step 4. In the inner loop, in Step 10 these parameters are respectively increased and decreased by factors of  $\Gamma$ and $\gamma$ after performing the Newton step. In Step 11, the accuracy of the new dual iterate $\bbmu_{k+1}$ is checked to be within the intended accuracy $\hat{V}$. Note that, while the sub-optimality cannot be checked directly without knowledge of $\tdR_{k+1}^*$, it can be checked indirectly by checking the norm of the gradient $\| \nabla \tdR_{k+1}(\bbmu_{k+1})\| < (\sqrt{2\alpha})\hat{V}$ from the strong convexity property in \eqref{eq_hessian_bounds}. If the condition in Step 11 is satisfied, the parameters $n$ and $M$ require no further modification. Otherwise, they are further modified until $\mu_{k+1}$ is within the target accuracy which in turn may imply that the conditions in \eqref{cond_1} and \eqref{cond_2} are satisfied. Note that the backtracking rates $\gamma, \Gamma$ are standard parameters used in the definition of a backtracking algorithm and effectively tradeoff the speed of the backtracking search vs. its thoroughness or accuracy. Generally speaking, values closer to 1 will result in a slower, more careful backtracking search while values of $\gamma$ and $\Gamma$ that are, respectively, smaller and larger will result in a faster, more aggressive search. Tuning of these parameters should thus reflect the desired tradeoff. With this practical considerations in mind, we proceed by simulating a wireless control learning problem using the proposed use of Newton's method on the ERM relaxation.

 
 \section{Simulation Results}\label{sec_numerical_results}
 
\begin{figure*}[t]
\centering
\includegraphics[width=0.31\linewidth,height=0.23\linewidth]{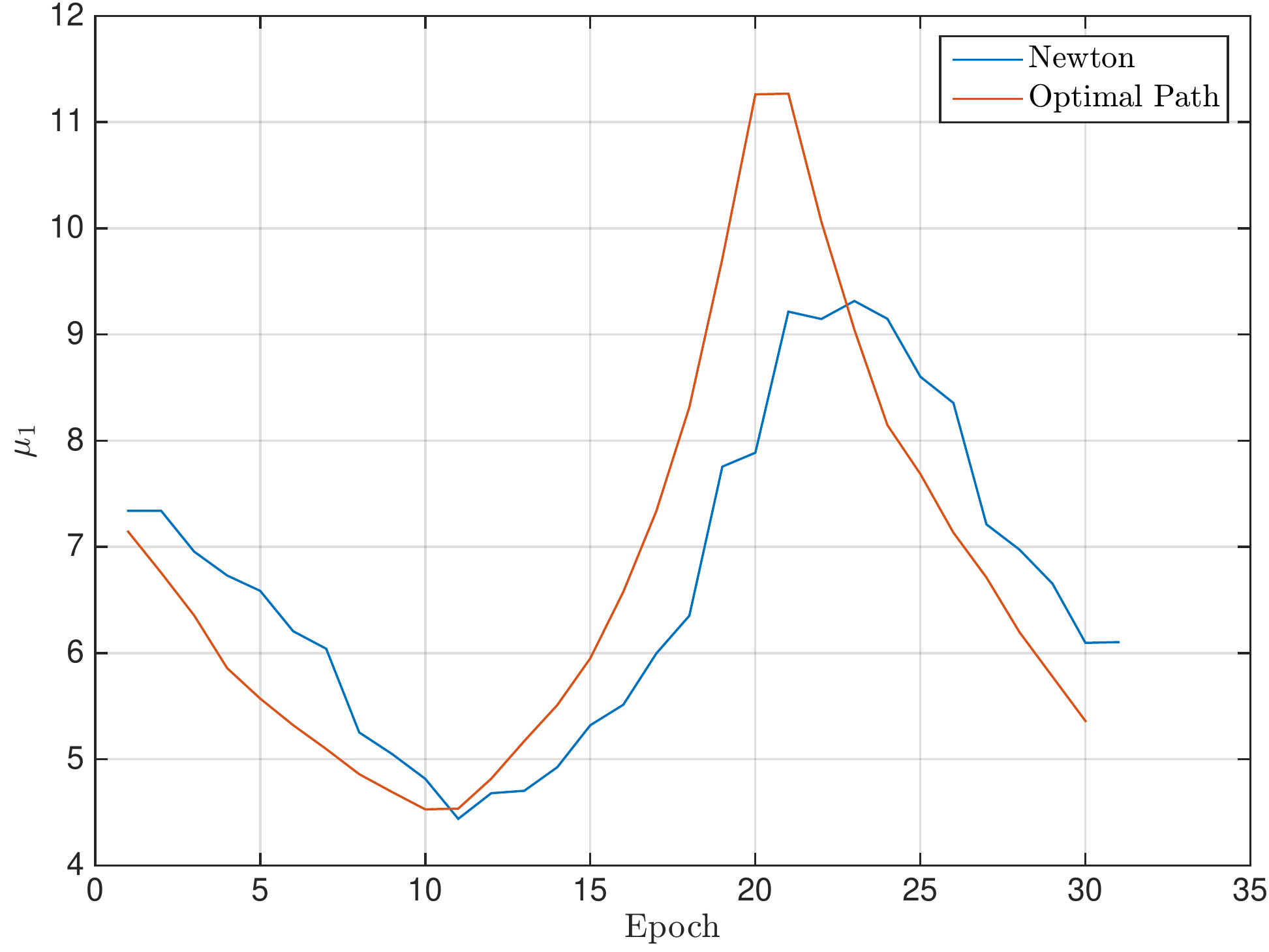}\quad
\includegraphics[width=0.31\linewidth,height=0.23\linewidth]{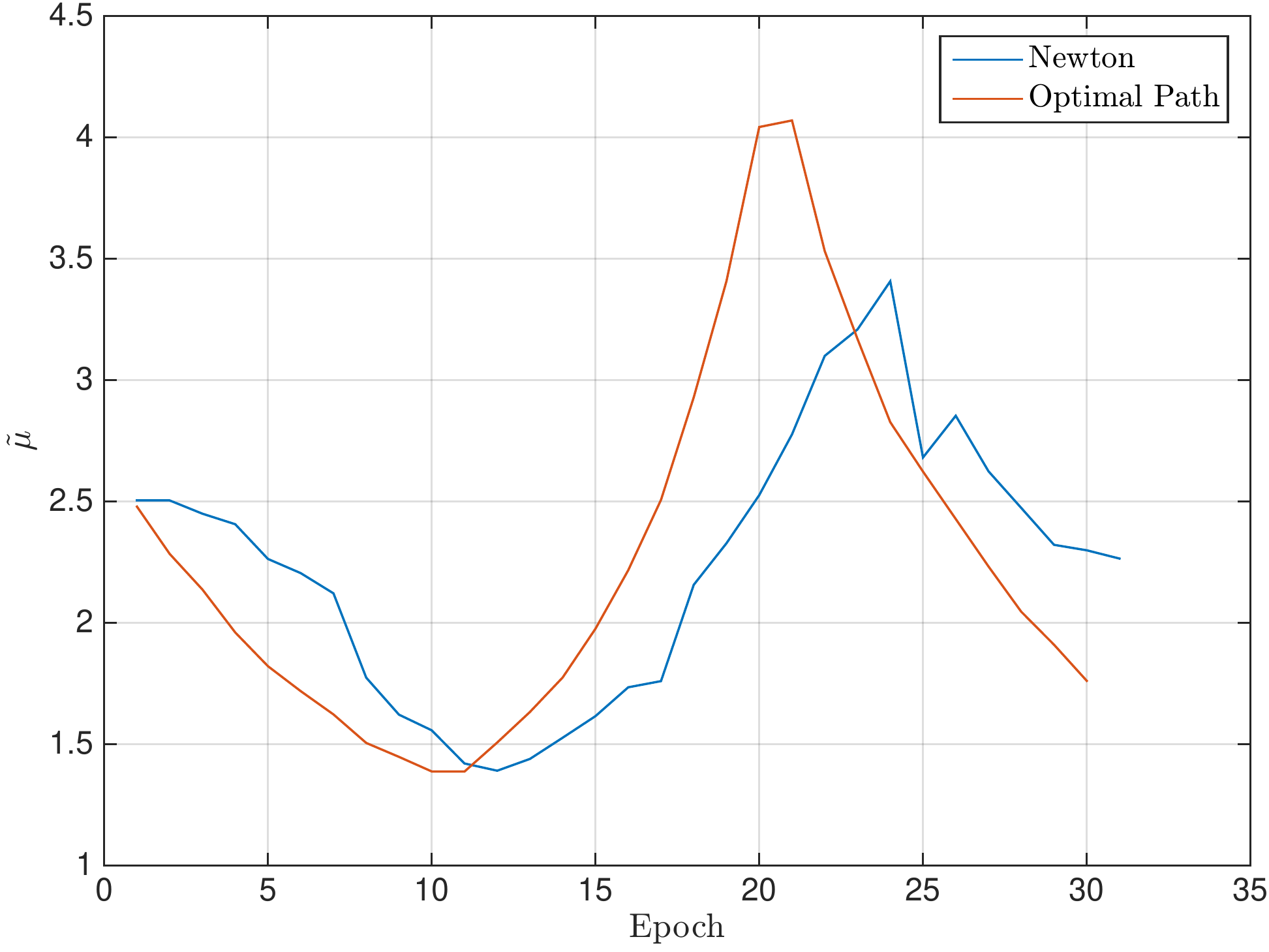} \quad
\includegraphics[width=0.31\linewidth,height=0.23\linewidth]{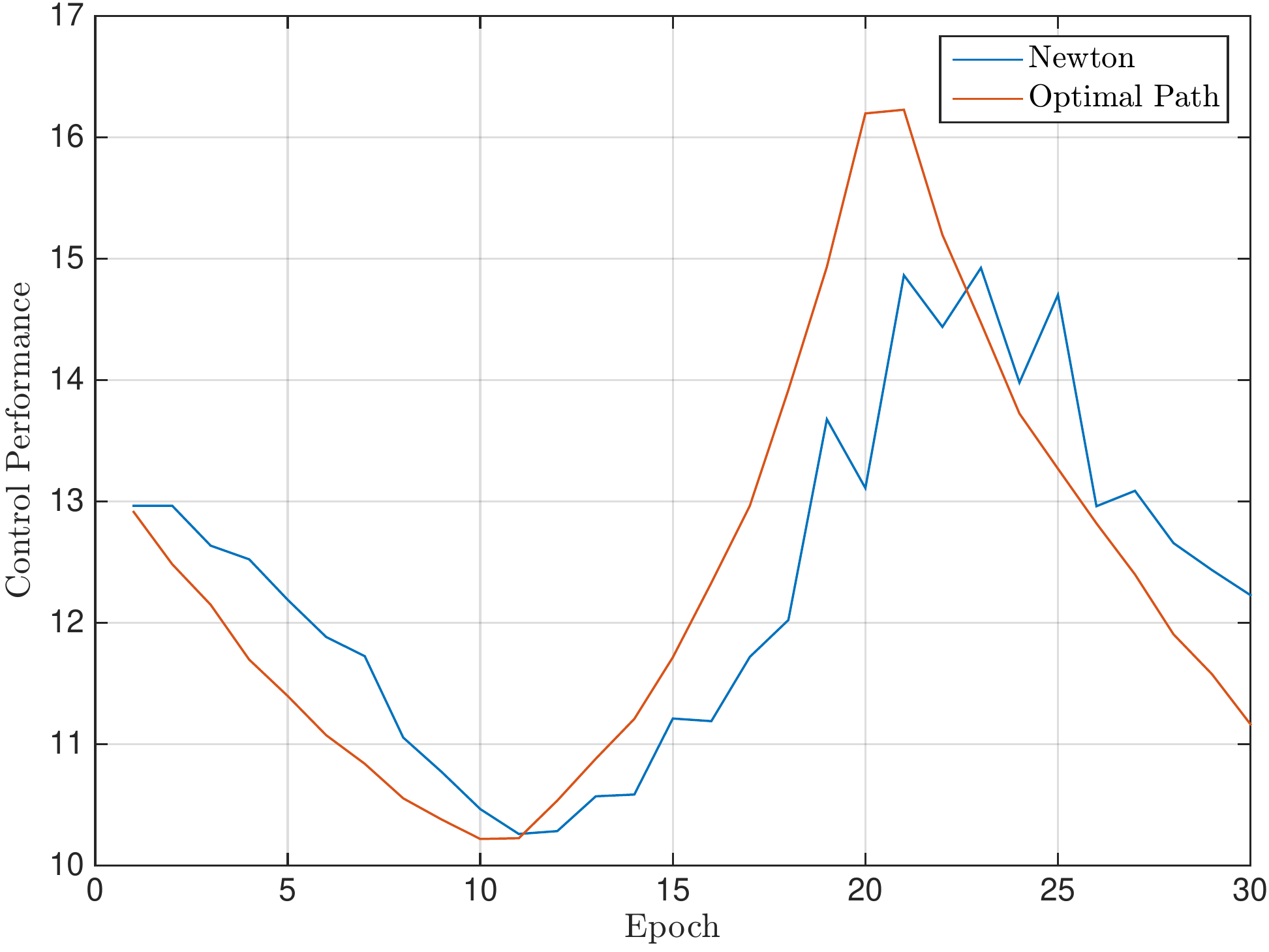}
\caption{Convergence paths of optimal values vs. values generated by the Newton learning method for time-varying $\ccalH_k$ for dual variables (left) $\mu^1$, (center) $\tdmu$, and (right) control performance $\sum J^i(y^i)$. Newton's method is able to find an approximately optimal value for the dual variables and respective control performance at each iteration. }\label{fig_sim_results}
\end{figure*}

 
We simulate the performance of our second order learning method on a simple WCP. Consider the 1-dimensional switched dynamical system in \eqref{example1} governed by the transition constants $A_o$ and $A_c$ for $m=4$ systems/states. The control performance for the $i$th agent $J^i(y^i)$ measures the mean square error performance and is now given by the expression in \eqref{eq_j_ex1}.
The open and closed loop control gains for each agents are chosen between $[1.1,1.5]$ and $[0,0.8]$, respectively. The probability of successful transmission for agent $i$ is modeled as a negative exponential function of both the power and channel state, $q(h^i,p^i(h^i) := 1 - e^{ -h^i p^i(h^i)}$, while channel states at epoch $k$ are drawn from an exponential distribution with mean $u_{k}$. The channel varies over time by the mean $u_k$ changing for different times. We draw $n=200$ samples and store a window of the previous $M=5$ distributions for a total of $N=1000$ samples at each epoch. As we assume the that channel statistics vary only vary across time \emph{epochs}, but stay constant within a single epoch, we may consider it reasonable to collect 200 channel samples within an epoch.

To demonstrate the ability of Newton's method to instantaneously learn an approximately optimal power allocation as the channel distribution varies over time, we perform Algorithm \ref{alg:AdaNewton} over the ERM problem in \eqref{eq_erm_problem} with the defined control performance $J(\cdot)$, transmission probabilities $q(\cdot)$ and channel distributions $\ccalH_k$. In Figure \ref{fig_sim_results} we show the path of Newton's method at each time $k$ for the dual variables $\mu^1_k$, $\tdmu_k$, and the control performance $\sum_{i=1}^m J^i(y^i_k)$. The red line of each figure plots the optimal values for the current distribution parameter $u_k$ as it changes with $k$. These values are obtained by solving the optimization problem at each epoch offline a priori. The blue line, alternatively, plots the values generated by Newton's method for each epoch $k$ in an online manner. The channel evolves at each iteration by a fixed rate $u_{k+1} = u_{k} \pm r$ for some rate $r$. Observe that within some small error Newton's method is indeed able to quickly and approximately find each new solution as the channel varies over time.

\begin{figure}
\centering
\includegraphics[width=\linewidth,height=0.3\textheight]{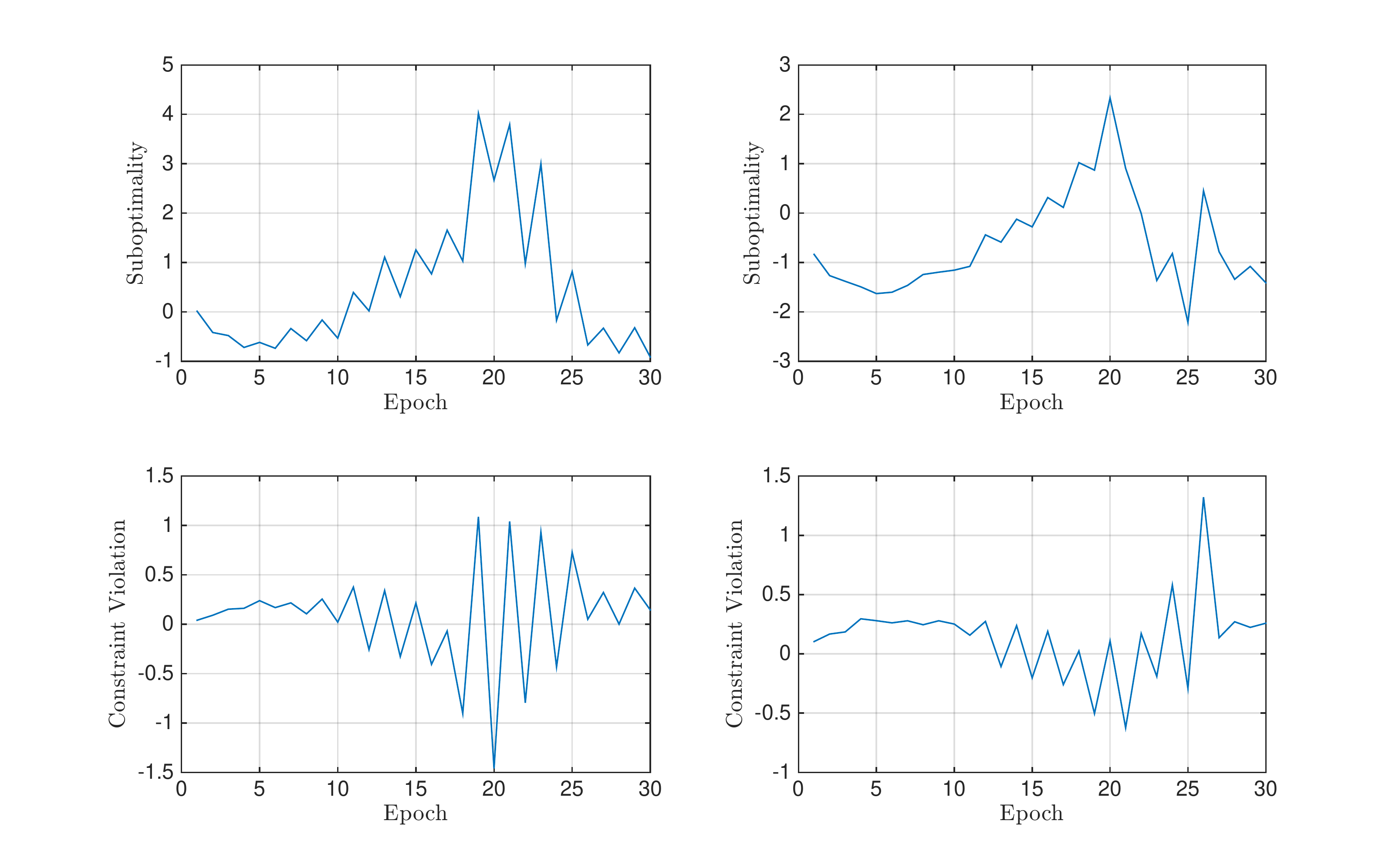}\quad
\caption{Comparison of suboptimality (top) and constraint violation (bottom) for the case of $\hat{V} = 0.01$ (left) and $\hat{V} = 0.03$ (right). Although the right-hand figures strive for less accuracy, they perform better because Newton's method can adapt to the intended accuracy more easily with single iterations. }\label{fig_stat_compare}
\end{figure}

To compare the effect of selecting different choices of accuracy $\hat{V}$ numerically, we present in Figure \ref{fig_stat_compare} the simulation performance of two representative cases with respective accuracies of  $\hat{V} = 0.01$ (left) and $\hat{V} = 0.03$ (right). In the top figures, we show the suboptimality relative to the optimal control performance and show on the bottom figures the resulting constraint violation (where a positive value reflect violation) over a set of time epochs where the channel varies. Here, we see an interesting case that highlights the need of proper selection or estimation of $\hat{V}$. Although the left hand figures strive for a better accuracy, the performance is better on the right hand figures. This is due to the fact that single iterations of Newton's method cannot reach accuracies of 0.01, resulting in a more suboptimal trajectory of resource allocation policies. On the other hand, the more moderate goal of 0.03 allows for the learning method to reach intended accurate goals with each step of Newton's method as the channel varies.

\begin{figure}
\centering
\includegraphics[width=\linewidth,height=0.3\textheight]{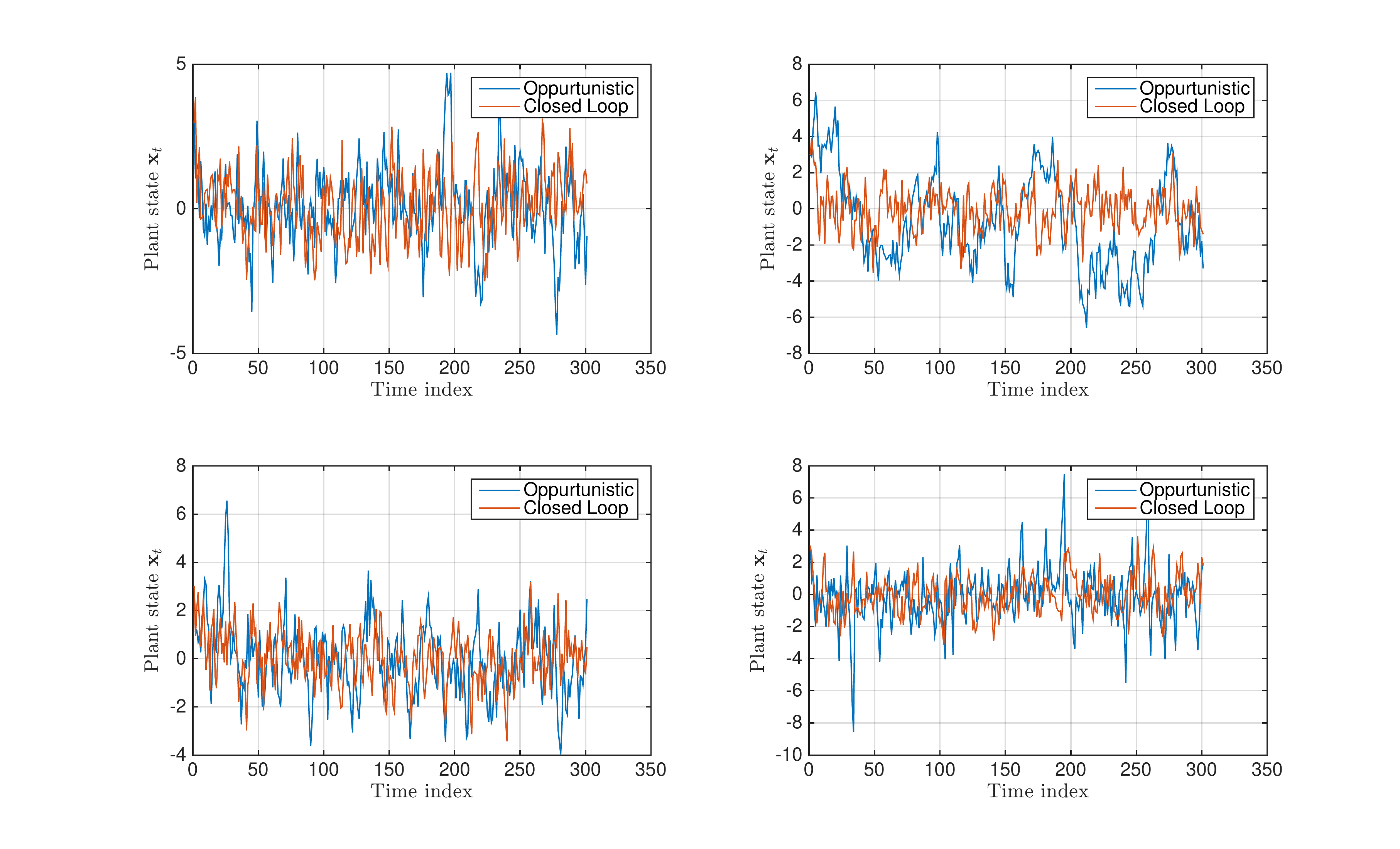}\quad
\caption{Dynamic evolution of each of the 4 state variables over the time-varying channel. The blue curve shows the opportunistic power allocation policy found with Newton's method while the red curve shows the evolution assuming the loop can always be closed. }\label{fig_sim_process}
\end{figure}

Using the dual parameters found by Newton's method, we simulate the resulting dynamical system. The dual parameters are used to determine the power allocation policy, which is used to determine transmission probabilities given current channel conditions. In Figure \ref{fig_sim_process} we show the resulting state evolution of $x^i_t$ for each of the 4 state variables. The blue curve shows the process using the opportunistic transmission policy from Newton's method, while the red curve shows the process when the loop is always closed, i.e. no packet drops. Here, we observe that while there are some instances when the state variable grows large when the system is in open loop, overall the system remains stable over time.

 \section{Conclusion}
In this paper we considered the wireless control system over a non-stationary wireless channel. The problem of maximizing a control utility subject to resource constraints can be formulated as a stochastic optimization problem in the dual domain. Because the wireless channel is random and time-varying, channel samples must be taken, resulting in a relaxed empirical risk minimization (ERM) problem. Standard ERM techniques do not suffice in the wireless setting because the channel is constantly changing. We propose the use of Newton's method, whose local quadratic convergence property allows us to continuously learn and adapt our optimal power allocation policies to changes in the channel distribution. We derive specific conditions on achieving instantaneous convergence to an approximate solution and characterize the suboptimality and stability in the wireless control problem (WCP). We additionally provide numerical simulations that demonstrate the use of Newton's method to learn and track optimal power allocations over a time varying channel. While this paper considers only resource allocation on contention-free links, consider the scheduling problem on a shared channel with non-stationary distributions remains an area of future work.

 
 \section*{Appendix: Proof of Lemma \ref{main_lemma}}\label{sec_main_lemma}
 
 We start with the definition of the Newton decrement at time $k+1$. We can add and subtract  $ \nabla \tdR_{k}(\bbmu_k)$ and upper bound using the triangle inequality as
\begin{align}\label{eq_1}
&\lambda_{k+1} (\bbmu_k) =  \|\bbH_{k+1}^{-1/2}  \nabla \tdR_{k+1}(\bbmu) \| = \| \nabla \tdR_{k+1}(\bbmu_k) \|_{\bbH_{k+1}^{-1}} \nonumber \\
 &\leq \| \nabla \tdR_{k}(\bbmu_k) \|_{\bbH_{k+1}^{-1}} + \| \nabla \tdR_{k+1}(\bbmu_k) - \nabla \tdR_{k} (\bbmu_k) \|_{\bbH_{k+1}^{-1}}.
\end{align}
First, we will upper bound the second term in \eqref{eq_1}. By adding and subtracting the expected losses $\nabla L_{k} (\bbmu_k)$ and $\nabla L_{k+1}(\bbmu_k)$ and using the triangle inequality to obtain
\begin{align}
& \| \nabla \tdR_{k+1}(\bbmu_k) - \nabla \tdR_{k} (\bbmu_k) \|
 \leq  \| \nabla \hat{L}_{k+1}(\bbmu_k) -\nabla L_{k+1} (\bbmu_k) \|\nonumber\\
 &\quad+  \|\nabla L_{k} (\bbmu_k) -  \nabla \hat{L}_k(\bbmu_k)  \| +  \| \nabla L_{k+1} (\bbmu_k) - \nabla L_{k} (\bbmu_k) \|. \nonumber
\end{align}
The first two terms in the above sum are bounded by $\tdV_N^{1/2}$ per \eqref{eqn_loss_minus_erm_2}, while the third term is the difference of two consecutive loss functions and is therefore bounded by $\bar{D_k}$ from \eqref{eqn_nonstat_bound_grad}. The norm weight $\bbH_{k+1}^{-1}$ additionally provides a bound of $\alpha \hat{V}$ as the strong convexity constant of $\tdR_{k+1}$ providing an upper bound on the norm of Hessian inverse as in \eqref{eq_hessian_bounds}. Combining these, we obtain
\begin{align}\label{eq_11}
 \| \nabla \tdR_{k+1}(\bbmu_k) - \nabla \tdR_{k} (\bbmu_k) \|_{\bbH_{k+1}^{-1}}   &\leq \frac{2 \tdV_{N}^{1/2} + \bar{D_k}}{(\alpha \hat{V})^{1/2}}.
\end{align}
We now can bound the first term in \eqref{eq_1} using the Lipschitz continuity of the gradient $\Delta + c\hat{V}$, i.e.
\begin{align}\label{q12}
\| \nabla \tdR_{k}(\bbmu_k) \|_{\bbH_{k+1}^{-1}}  \leq \left( \frac{2(\Delta+c\hat{V})\|\bbmu_{k}-\tbmu_k^*\|}{\alpha\hat{V}}\right)^{1/2}
\end{align}
Recall that $\bbmu_k$ is given to be a $\hat{V}$-accurate minimizer of $\tdR_k$. The difference $\|\bbmu_{k}-\tbmu_k^*\|$ can subsequently be bounded with $\hat{V}$, resulting in the final bound for the first term
\begin{align}\label{q1}
\| \nabla \tdR_{k}(\bbmu_k) \|_{\bbH_{k+1}^{-1}}  \leq \left( \frac{2(\Delta+c\hat{V})\hat{V}}{\alpha\hat{V}}\right)^{1/2}
\end{align}
To be in the quadratic convergence region, i.e. $\lambda_{k+1}(\bbmu_k) < 1/4$, follows by summing \eqref{eq_11} and \eqref{q1} as in \eqref{prop_result}.

\section*{Appendix: Proof of Lemma \ref{lemma_sub_m}}\label{sec_lemma_sub_m}
To prove this result, we start by expanding the term $\tdR_{k+1}(\bbmu_k) - \tdR_{k+1}^*$. By adding and subtracting $\tdR_{k}(\bbmu_{k})$, $\tdR_{k}^*$, and $\tdR_{k}(\bbmu_{k+1}^*)$, we obtain
\begin{align}\label{proof_eq_2_11}
 \tdR_{k+1}(\bbmu_{k}) - \tdR_{k+1}^*
 &=  \tdR_{k+1}(\bbmu_{k}) -  \tdR_{k}(\bbmu_{k})\\
&\qquad +  \tdR_{k}(\bbmu_{k}) - \tdR_{k}^* \nonumber\\
 &\qquad+ \tdR_{k}^* -\tdR_{k}(\bbmu_{k+1}^*) \nonumber\\
 &\qquad+\tdR_{k}(\bbmu_{k+1}^*) - \tdR_{k+1}^*. \nonumber
\end{align}
We now individually bound each of the four differences in \eqref{proof_eq_2_11}. Firstly, the difference $\tdR_{k+1}(\bbmu_{k}) -  \tdR_{k}(\bbmu_{k})$ becomes 
\begin{align}\label{proof_eq_2_22}
\tdR_{k+1}(\bbmu_{k}) -  \tdR_{k}(\bbmu_{k})
&=\hat{L}_{k+1}(\bbmu_{k})-\hat{L}_{k}(\bbmu_{k}),
\end{align}
Using the same reasoning as in \eqref{eq_11} with the functional statistical accuracy bound in place of the bound for gradients in \eqref{eqn_loss_minus_erm_2} and using \eqref{eqn_nonstat_bound} in place of \eqref{eqn_nonstat_bound_grad}, we obtain the equivalent bound
\begin{equation}\label{proof_eq_2_44}
\tdR_{k+1}(\bbmu_{k}) -  \tdR_{k}(\bbmu_{k})  \leq   2\tdV_N + D_k.
\end{equation}
For the second term in \eqref{proof_eq_2_11}, we again use the fact that $\bbmu_{k}$ as an $\hat{V}$-optimal solution for the sub-optimality $ \tdR_{k}(\bbmu_{k}) - \tdR_{k}^* $ to bound with the statistical accuracy as 
\begin{equation}\label{proof_eq_2_55}
 \tdR_{k}(\bbmu_{k}) - \tdR_{k}^*   \leq  \hat{V}.
\end{equation}
We proceed with bounding the third term in \eqref{proof_eq_2_11}. Based on the definition of $\bbmu_{k}^*$ as the optimal solution of the loss $\tdR_{k}$, the the difference $ \tdR_{k}^* -\tdR_{k}(\bbmu_{k+1}^*) $ is always negative, i.e., 
\begin{equation}\label{proof_eq_2_66}
\tdR_{k}^* -\tdR_{k}(\bbmu_{k+1}^*)    \leq  0.
\end{equation}
For the fourth term in \eqref{proof_eq_2_11}, we use the triangle inequality to bound the difference $\tdR_{k}(\bbmu_{k+1}^*) - \tdR_{k+1}^*$ in \eqref{proof_eq_2_11} as
\begin{align}\label{proof_eq_2_77}
\tdR_{k}(\bbmu_{k+1}^*) - \tdR_{k+1}^* 
& =  \hat{L}_{k}(\bbmu_{k+1}^*)-\hat{L}_{k+1}(\bbmu_{k+1}^*) \nonumber\\
& \leq  2\tdV_N + D_k.
\end{align}
Observe that \eqref{proof_eq_2_77} uses the same reasoning as \eqref{proof_eq_2_44}. Replacing the differences in \eqref{proof_eq_2_11} by the upper bounds in \eqref{proof_eq_2_44}-\eqref{proof_eq_2_77},
\begin{equation}\label{proof_eq_2_88}
 \tdR_{k+1}(\bbmu_{k}) - \tdR_{k+1}^* \leq 4\tdV_N + \hat{V} + 2D_k
\quad \text{w.h.p.}
\end{equation}

\section*{Appendix: Proof of Lemma \ref{quadratic_convg_lemma}}\label{sec_quadratic_convg_lemma}
The proof for this result follows from \cite[Proposition 4]{mokhtari2016adaptive}, which we repeat here for completeness. We proceed by bounding the difference $\tdR_{k+1}(\bbmu)-\tdR_{k+1}^*$ in terms of the Newton decrement parameter $\lambda_{k+1} (\bbmu)$. We first use the result in \cite[Theorem 4.1.11]{nesterov2013introductory}, showing that
 \begin{align}\label{proof_of_prop_10}
		  \lambda_{k+1} (\bbmu) - \ln\left(1+\lambda_{k+1} (\bbmu)\right) &\leq  
	 \tdR_{k+1}(\bbmu)-\tdR_{k+1}^*   \\
	 &\leq
	  - \lambda_{k+1} (\bbmu) - \ln\left(1-\lambda_{k+1} (\bbmu)\right). \nonumber
\end{align}
We can use the Taylor's expansion of $\ln(1+a)$ for $a=\lambda_{k+1} (\bbmu)$ to show that $\lambda_{k+1} (\bbmu) - \ln\left(1+\lambda_{k+1} (\bbmu)\right)$ is bounded below by $(1/2)\lambda_{k+1} (\bbmu)^2-(1/3)\lambda_{k+1} (\bbmu)^3$ for $0<\lambda_{k+1} (\bbmu)<1/4$ . Likewise, we have that $(1/6)\lambda_{k+1} (\bbmu)^2\leq (1/2)\lambda_{k+1} (\bbmu)^2-(1/3)\lambda_{k+1} (\bbmu)^3$ and subsequently $\lambda_{k+1} (\bbmu) - \ln\left(1+\lambda_{k+1} (\bbmu)\right)$ is bounded below by $(1/6)\lambda^2$. We again use Taylor's expansion of $\ln(1-a)$ for  $a=\lambda_{k+1} (\bbmu)$ to show that $- \lambda_{k+1} (\bbmu) - \ln\left(1-\lambda_{k+1} (\bbmu)\right) $ is bounded above by $\lambda_{k+1} (\bbmu)^2$ for $\lambda_{k+1} (\bbmu)<1/4$; see e.g., \cite[Ch. 9]{boyd2004convex}. Considering these bounds and the inequalities in \eqref{proof_of_prop_10} we obtain that
 \begin{equation}\label{proof_of_prop_20}
	\frac{1}{6}\lambda_{k+1} (\bbmu)^2\leq
	 \tdR_{k+1}(\bbmu)-\tdR_{k+1}^*  \leq
	 \lambda_{k+1} (\bbmu)^2.
\end{equation}

Because we assume that $\lambda_{k+1} (\bbmu_k)\leq 1/4$, the quadratic convergence rate of Newton's method for self-concordant functions \cite{boyd2004convex} implies that the Newton decrement has a quadratic convergence and we can write 
 \begin{equation}\label{proof_of_prop_30}
\lambda_{k+1} (\bbmu_{k+1}) \leq 2 \lambda_{k+1} (\bbmu_k)^2. 
\end{equation}
We combine the results in \eqref{proof_of_prop_20} and \eqref{proof_of_prop_30} to show that the optimality error $\tdR_{k+1}(\bbmu_{k+1})-\tdR_{k+1}^*$ has an upper bound which is proportional to $(\tdR_{k+1}(\mu_k)-\tdR_{k+1}^* )^2$. In particular, we can write $ \tdR_{k+1}(\bbmu_{k+1})-\tdR_{k+1}^* \leq \lambda_{k+1} (\bbmu_{k+1})^2 $ based on the second inequality in \eqref{proof_of_prop_20}. This observation in conjunction with the result in \eqref{proof_of_prop_30} implies that
 \begin{align}\label{proof_of_prop_40}
 \tdR_{k+1}(\bbmu_{k+1})-\tdR_{k+1}^* \leq 4 \lambda_{k+1} (\bbmu_k)^4.	
\end{align}
The first inequality in \eqref{proof_of_prop_20} implies that $\lambda_{k+1} (\bbmu_k)^4\leq 36 (\tdR_{k+1}(\bbmu_k)-\tdR_{k+1}^*)^2 $. Thus, we can substitute $\lambda_{k+1} (\bbmu_k)^4$ in \eqref{proof_of_prop_40} by $36 (\tdR_{k+1}(\bbmu_k)-\tdR_{k+1}^*)^2 $ to obtain the result in \eqref{prop_result_2}.

\section*{Appendix: Proof of Proposition \ref{feasible_error_prop}}\label{sec_feasible_error_prop}
We begin by bounding the gradient of the expected dual loss $L(\bbmu_k)$ at the $k$th dual iterate $\bbmu_k$ by using Lipschitz continuity, i.e.
\begin{align}
\| \nabla L_k (\bbmu_k) \|^2 \leq 2\Delta (L_k(\bbmu_k) - L_k^*).
\end{align}
We expand the sub-optimality $L(\bbmu_k) - L^*$ by adding and subtracting terms as follows
\begin{align}\label{eq_pp_1}
\frac{1}{2\Delta}\| \nabla L_k (\bbmu_k) \|^2 &\leq L_k(\bbmu_k) - \tdL_k(\bbmu_k) + \tdL_k(\bbmu_k)  \\&\quad - \tdR_k(\bbmu_k) + \tdR_k(\bbmu_k) - \tdR^*_k + \tdR^*_k - L_k^*, \nonumber
\end{align}
where we recall the notation $\tdR_k^* := \tdR_k(\tbmu_k^*)$. We now proceed by bounding each successive pair of terms in \eqref{eq_pp_1}. The first difference $L_k(\bbmu_k) - \tdL_k(\bbmu_k)$ comes from the sampling and is thus bounded by the statistical accuracy $\tdV_N$. The second difference $\tdL_k(\bbmu_k)  - \tdR_k(\bbmu_k)$ can be bounded by the regularizers as
\begin{align}\label{eq_pp}
\tdL_k(\bbmu_k)  - \tdR_k(\bbmu_k) \leq \beta \hat{V} \mathbf{1}^T \log_{\eps} (\bbmu_k) - \frac{\alpha \hat{V}}{2}\| \bbmu_k\|^2.
\end{align}
The second term on the right hand side of \eqref{eq_pp} is negative and can be ignored. Because the dual variable $\| \bbmu_k\|$ was upper bounded in Corollary \ref{corollary_primal_bound}, we can place a finite bound on $\mathbf{1}^T \log_{\eps} (\bbmu_k) \leq \kappa$ and then bound the term $\beta \hat{V} \mathbf{1}^T \log_{\eps} (\bbmu_k) \leq \beta \hat{V} \kappa$.  The third difference $\tdR_k(\bbmu_k) - \tdR^*_k$ is bounded by the suboptimality $\hat{V}$ from the main result in \eqref{imp_result} and the fourth difference $\tdR^*_k - L_k^*$ can be bounded by $\rho \hat{V}$ from \eqref{eq_prop_diff}. We can therefore bound the gradient of the dual loss as
\begin{align} \label{eq_pp2}
\| \nabla L_k (\bbmu_k) \|^2 \leq 2\Delta(\tdV_N + C \hat{V}),
\end{align}
where $C := 1 + \rho + \beta m \log{\kappa}$. From here, consider that the norm of the dual gradient $\| \nabla L_k (\bbmu_k) \|^2$ is the sum of squares of each constraint violation in \eqref{eq_power_problem_y}, i.e.,
\begin{align}
&\left(\sum_{i=1}^m \E_{h_k^i} (p^i(h)) - p_{\max}\right)^2 + \sum_{i=1}^m \left(y^i -  \E_{h_k^i} \left\{q(h,p^i(h))\right\} \right)^2 \nonumber\\&\quad \leq 2\Delta (\tdV_N  + C \hat{V}).
\end{align}
The results in \eqref{eq_j1} and \eqref{eq_j2} then follow from here.

\section*{Appendix: Proof of Theorem \ref{cor_primal_error}}\label{sec_cor_primal_error}
Consider that, using the definitions of the primal maximizers $\bbp(\bbh,\bbmu_k)$ and $\bby(\bbmu_k)$ at a dual point $\bbmu_k$, we can write the dual function as
\begin{align}
L(\bbmu_k) = J(\bby(\bbmu_k)) + \bbmu_k^T\left(\E_{\bbh} \check{\bbq}(\bbp(\bbh,\bbmu_k)) - \check{\bby}(\bbmu_k)\right).
\end{align}
Likewise, we know from strong duality that the optimal dual values $L_k^*$ is equivalent to the optimal primal objective value $J(\bby_k^*)$. Therefore, we can write the suboptimality of dual functions as
\begin{align} \label{eq0001}
L(\bbmu_k) - L_k^* &= J(\bby(\bbmu_k)) - J(\bby_k^*) \\ &\qquad + \bbmu_k^T\left(\E_{\bbh} \check{\bbq}(\bbp(\bbh,\bbmu_k)) - \check{\bby}(\bbmu_k)\right). \nonumber
\end{align}
Using the bound on dual suboptimality that comes from combining strong convexity and the gradient bound in \eqref{eq_pp2}, we can upper bound \eqref{eq0001} as
\begin{align} \label{eq0002}
(1+C)\Delta/\alpha &\geq J(\bby(\bbmu_k)) - J(\bby_k^*) \\ &\qquad + \bbmu_k^T\left(\E_{\bbh} \check{\bbq}(\bbp(\bbh,\bbmu_k)) - \check{\bby}(\bbmu_k)\right). \nonumber
\end{align}
We can lower bound the right hand side of \eqref{eq0002} by taking the negative of the absolute value of the final term. Rearranging terms we obtain
\begin{align} \label{eq0003}
(1+C)\Delta/\alpha + &|\bbmu_k^T\left(\E_{\bbh} \check{\bbq}(\bbp(\bbh,\bbmu_k)) - \check{\bby}(\bbmu_k)\right)| \\& \geq J(\bby(\bbmu_k)) - J(\bby_k^*).\nonumber
\end{align}
From here, we can upper bound the second term on the left hand side using the Cauchy-Schwartz inequality. The norm $\|\bbmu_k\|$ is bounded by $\sqrt{2/\alpha} + \hat{\ccalK}$ from Corollary \ref{corollary_primal_bound} while the norm $\|\E_{\bbh} \check{\bbq}(\bbp(\bbh,\bbmu_k)) - \check{\bby}(\bbmu_k)\|$ is bounded by $2\Delta (1+C) \hat{V}$ from \eqref{eq_j2}. This provides us the final result as
\begin{align} \label{eq0004}
(1+C)\Delta\left(\frac{1}{\alpha} + 2 \hat{V} (\sqrt{2/\alpha} + \hat{\ccalK})\right) \geq J(\bby(\bbmu_k)) - J(\bby_k^*).
\end{align}
%

\urlstyle{same}
\bibliographystyle{IEEEtran}
\bibliography{../wireless_control}

\end{document}